\documentclass[11pt, twosides]{amsart}
\usepackage{defs}
\usepackage{amsaddr}

\title[Stabilization under round robin scheduling of controllers]{Stabilization under round robin scheduling of control inputs in nonlinear systems}
\author[C. Maheshwari, S. Srikant, and D. Chatterjee]{Chinmay Maheshwari}
\address{EECS\\ University of California Berkeley\\ Berkeley 94720, USA}
\author{Sukumar Srikant and Debasish Chatterjee}
\address{Systems and Control Engineering\\ IIT Bombay, Powai\\ Mumbai 400076, India\\ \url{http://www.sc.iitb.ac.in/~srikant}\\ \url{http://www.sc.iitb.ac.in/~chatterjee}}
\email{{chinmay\_maheshwari@berkeley.edu,\{srikant.sukumar,dchatter\}@iitb.ac.in}}

\begin{document}

\begin{abstract}
	We study stability of multivariable control-affine nonlinear systems under sparsification of feedback controllers. Sparsification in our context refers to the scheduling of the individual control inputs one at a time in rapid periodic sweeps over the set of control inputs,  which corresponds to round-robin scheduling. We prove that if a locally asymptotically stabilizing feedback controller is sparsified via the round-robin scheme and each control action is scaled appropriately, then the corresponding equilibrium of the resulting system is stabilized when the scheduling is sufficiently fast; under mild additional conditions, local asymptotic stabilization of the corresponding equilibrium can also be guaranteed. {Moreover, the basin of attraction for the equilibrium of scheduled system also remains same as the original system under sufficiently fast switching. } Our technical tools are derived from optimal control theory, and our results also contribute to the literature on the stability of switched systems in the fast switching regime. Illustrative numerical examples depicting several subtle features of our results are included.
\end{abstract}


\maketitle

%%%%%%%%%%%%%%%%%%%%%%%%%%%%%%%%%%%%%%%%%%%%%%%%%%%%%%%%%%%%%%%%%%%%%%%%%%%%%%%%
%%%%%%%%%%%%%%%%%%%%%%%%%%%%%%%%%%%%%%%%%%%%%%%%%%%%%%%%%%%%%%%%%%%%%%%%%%%%%%%%
\section{Introduction}
\label{sec: Intro}
%%%%%%%%%%%%%%%%%%%%%%%%%%%%%%%%%%%%%%%%%%%%%%%%%%%%%%%%%%%%%%%%%%%%%%%%%%%%%%%%
For positive integers \(d\) and \(m\), consider a control-affine nonlinear system with $m$ control inputs 
\begin{align}\label{eq: cfa}
	\dot{x}(t) = f(x(t)) + \sum_{i=1}^{m}g_i(x(t))u_i(t),
\end{align}
where $f: \R^d \rightarrow \R^d$ is a smooth drift vector field, $g_i : \R^d \rightarrow \R^d$ are smooth control vector fields for each \(i\), and \(x(t) \in \R^d \) is the vector of states and $u_i(t) \in \R$ is the \(i\)-th control action at time \(t\). { A major chunk of the literature on control theory has the underlying assumption that } all the \(m\) control inputs are at the disposal of the control designer at every instant of time.
{ However, this assumption is restrictive for large scale networked control systems \cite{walsh2001scheduling,gorges2007optimal,Nesic2004Automatica,Nesic2004TAC,Heemels2019periodic}  where the actuation and sensing resources are shared between various sub-systems in order to reduce the cost and complexity. In such situations, often the actuators and sensors are scheduled such that not all of them are continously updated. 
}
%In large scale networked control systems, typically, the actuation and sensing capabilities are resources shared between various sub-systems, and therein such scheduling is motivated by the need to reduce cost and complexity.
Moreover, simultaneously commanding $m$ control inputs at every instant of time for high-dimensional systems ($d, m \gg 1$) is becoming a challenging task \cite[p.\ 111]{JurdQuinn} today, and such systems routinely arise in the context of multi-agent systems \cite{JurdQuinn}, traffic networks \cite{chen2010review}, smart grids \cite{SG1,claes2011decentralized}, etc.

%{\color{cyan}\sout{ A common feature in such systems is the fact that the active engagement of all the control inputs simultaneously leads to undesirable consequences such as increased communication costs, higher maintenance requirements, wireless signal interference, lack of flexibility, etc. Naturally, }}
On the one hand, a reduction in the number of control inputs that are simultaneously active at every time instant, is a desirable feature. On the other hand, reducing the active control inputs without proper care may be disastrous, leading to the loss of key system-theoretic properties such as controllability. Indeed, determining the fewest actuators for linear systems while retaining their controllability is a hard problem from a computational complexity standpoint; see \cite{polyak2013lmi} and the references therein.  

In this article we take advantage of sparsity to reduce the number of control actions that are simultaneously active. We think of each (scalar) control variable in a multidimensional controller as a separate entity (see Remark \ref{rem: allChannels} for a generalization), and sparsify the control actions by
\begin{itemize}[leftmargin=*,label=\(\circ\)]
	\item switching, at each time instant, \textbf{only} one of the inputs to the `on' or active state, while the rest are forced to attain the zero value (their `off' or inactive states), and 
	\item cycling rapidly over the set of inputs in a round-robin fashion.
\end{itemize}
The article at hand studies the effect of the preceding scheduling scheme in the context of stabilization of the nonlinear system \eqref{eq: cfa}.

To this end, first, we design an asymptotically stabilizing \emph{nominal} feedback controller for an equilibrium, say \(0\), of \eqref{eq: cfa} assuming that all the control inputs are always active. We do \emph{not} revisit the topic of synthesizing asymptotically stabilizing feedbacks for control-affine nonlinear systems since it has been extremely well studied in the literature, and instead assume that the engineer may pick a nominal controller off the shelf. Next, the following scheme is adopted:
\begin{enumerate}[leftmargin=*, label = (\Roman*), align=right, widest=ii]
	\item \label{intro: 1} we scale all the control values of the nominal controller by \(m\) (recall that \(m\) is the number of controllers), and
	\item \label{intro: 2} we pick a small number \(\tau >0\), select an arbitrary order of the \(m\) controllers (relabel them to be \((1,2,\dots, m)\) if necessary), and cyclically switch, every \(\tau \) units of time, one of the (scaled) controllers to `on' and the rest to `off' during the corresponding interval of length \(\tau \).
\end{enumerate}
The resulting cyclic and periodic multiplexing of the control actions leads to a particular kind of on-off behavior of a part of the closed-loop vector field, thereby producing a switched dynamical system in which \(\tau^{-1}\) plays the role of the \emph{switching frequency}. In the preceding setting, we shall establish that if the switching frequency is held fixed at a sufficiently high value (equivalently, if \(\tau\) is sufficiently small), then one can guarantee stability of the point \(0\). Moreover, if the switching is made faster with time in an appropriate way (to be made precise in \S\ref{ssec: AS}), then, in addition, asymptotic stability of \(0\) of the resulting switched system can be certified. We reiterate that our results depend on the mild requirement that closed-loop system has \(0\) to be a locally asymptotically stable equilibrium under the nominal controller. { Our final set of results, provided in \S\ref{ssec:BOA}, concerns precise estimates of the basin of attraction of asymptotically stable equilibrium points of the nominal system under round-robin multiplexing of the controllers.}

%\sout{The amplification aspect of our scheme discussed in \ref{intro: 1} above is a very crucial requirement for our results. Intuitively, since only one control input is active at any given instant, the amplification factor compensates for their inactive periods by depleting sufficient energy away from the system during its active period, thereby promoting stabilizing behaviour. We shall illustrate and elaborate more on the technical aspects of this scaling in \S \ref{sec: NumExp}. The amplification is, in spirit, related to the ideas in \cite{sukumar2011precision} where the authors designed a high-gain feedback controller to schedule a single actuator for 3-axis stabilization of a spacecraft. }

The emergence of a switched dynamical system from our scheduling scheme leads to new challenges for at least two reasons. One, it is well known that unplanned switching between systems may lead to instability \cite{LibBook} and even chaotic behaviour \cite{el2008chaotic}. Indeed, such undesirable situations may arise in our context if one is not careful, as will be illustrated by our numerical case studies presented in \S \ref{sec: NumExp}. Two, the fast switching regime has been only sparsely investigated in the literature, and there are practically no standard results to lift off the shelf. Our technical proofs proceed ab ovo, and our key tool is borrowed from optimal control in the form of the so-called chattering lemma \cite[Chapter 3]{ref:BerMed-13}; in our context, this lemma assures uniform convergence, as the switching frequency diverges to infinity, of the trajectories generated by the periodically scheduled controller to the trajectories generated by the nominally controlled system.

{The strategy of first designing a feedback controller being oblivious of scheduling and then proving certain stability properties of the resulting system (a strategy commonly referred to as \emph{emulation}), is not new. In fact, it is a heavily studied topic in the area of networked control systems; see, e.g., \cite{Nesic2004TAC,Nesic2004Automatica,walsh2001asymptotic,walsh2001scheduling,JurdQuinn,van2017switched,SchenatoZeroing,schenato2007foundations,YukselOptimalControl}. However, these studies differ from ours in at least one of the following ways: 
	\begin{enumerate}[leftmargin=*,label=\(\circ\)]
		\item sparsification is absent, i.e., all the controllers are always active \cite{Nesic2004TAC},
		%\item{\color{cyan} \sout{work with weaker stability notion than ours (\cite{Nesic2004TAC}), }}\todo[inline]{I feel that we can remove this point from here}
		\item impose restrictions on the system dynamics in order to incorporate high-level tools \cite{walsh2001asymptotic,walsh2001scheduling,JurdQuinn,wang2012input,Nesic2004Automatica,Nesic2004TAC,wang2015emulation,donkers2011stability},
		\item conduct emprical evidence of merit of strategy on very specific systems (e.g., for SCARA robots in \cite{van2017switched}).
		\item {study the effect of zero-input (or \emph{zeroing}) strategy on the linear quadratic performance under control packet loss \cite{SchenatoZeroing,YukselOptimalControl,schenato2007foundations}}
	\end{enumerate}
	Complementary to the article at hand is \cite{ref:BacMaz-10} where authors develop time and state dependent (typically slow) switching such that the scheduled system is stabilized \cite{ref:BacMaz-10}; we discuss more about this work in Remark \ref{rem: BacMaz}. 
}

{ The results of \cite{JurdQuinn} are more aligned with our work; they cater to nonlinear systems and provide asymptotic stability certificates for control-affine systems of the form \eqref{eq: cfa} under round-robin periodic scheduling of Jurdjevic-Quinn controllers \cite{jurdjevic1978controllability}. However, the drift vector fields are required to be dissipative in \cite{JurdQuinn}; as a result, the uncontrolled dynamics lose ``energy'' with time, which means that the system exhibits inherently ``stable'' behavior under zero control. In this article we remove such restrictions on the drift vector fields; in broad strokes, the only indispensable assumption employed here is that there are asymptotically stabilizing controllers available for our control-affine system \eqref{eq: cfa}.}

Our article unfolds as follows. We set up the premise of this article in \S\ref{ssec: Premise}. Our main results are presented in \S\ref{sec:MRS}, where we provide stability certificates under round-robin periodic scheduling of control inputs {and the corresponding effect on the basin of attraction of the equilibrium point}. In \S\ref{sec: NumExp} we present two numerical case studies that illustrate our main results presented \S\ref{sec:MRS}. A brief summary of the results and future research directions are provided in \S\ref{sec: Conc}. 

%%==============================================================================
\subsection*{Notations}
%%==============================================================================
Standard notations are employed here. We denote the set of non-negative integers by \(\nat \) and the set of positive integers by \(\N \). For a positive integer \(k \) we define  \(\setInt[k] \Let \set[]{1,2,\dots,k} \). The Lesbesgue measure of a measurable set \(S \subset \R \) is written as \(\abs{S}\). The interior and closure of a set \(A\subset\R^n\) are denoted by \(\inLenInt[A]\) and \(\closure[A] \), respectively. For a vector \(x \in \R^n\) we write \(\norm{x}\) for its Euclidean norm, and for \(r > 0\) we write the \emph{closed} Euclidean ball of radius \(r\) centered at \(x\) as \(\ball{r}{x} \Let \set[\big]{y \in \R^n \suchthat \norm{y-x} \leq c} \). For \(y \in \R \) we employ \(\left \lfloor y \right \rfloor \) to denote the greatest integer dominated by \(y\). If \(v \in \R^n \) is given, then \(\textsf{diag}(v) \) denotes a matrix \(D\) of dimension \(n \times n \), such that \(D_{ij} = v_i\delta_{ij}, i \in \setInt[n] \), with \(\delta_{ij} \) being the Kronecker delta.  

%%%%%%%%%%%%%%%%%%%%%%%%%%%%%%%%%%%%%%%%%%%%%%%%%%%%%%%%%%%%%%%%%%%%%%%%%%%%%%%%
\section{Main Results}
\label{sec: MainResultHeading}
\subsection{Stability Definitions}
%%%%%%%%%%%%%%%%%%%%%%%%%%%%%%%%%%%%%%%%%%%%%%%%%%%%%%%%%%%%%%%%%%%%%%%%%%%%%%%%
Let us start by setting up the basic definitions needed for our main results. The word \emph{measurability} will always stand for Lebesgue measurability in what follows.

For \(d, \filPPS\in\N\) and \(t_0\in\R\), consider the system 
\begin{align}\label{eq: firstSys}
	\dot{x}(t) = f(t,x(t);\filParam(t)), \quad\ x(\initT) = \bar{x}, \quad t \geq \initT,
\end{align}
where \(x(t) \in \R^d \) is the vector of states at time \(t\), and \(\filParam: \lcro{\initT, +\infty} \rightarrow  \R^\filPPS \) is a time-dependent map describing the evolution of a parameter. We admit only such maps \(\filParam\) that \(\lcro{\initT, +\infty}\times\R^d \ni (t,x)\mapsto f(t,x;\filParam(t)) \in \R^d\) is measurable in \(t\) and locally Lipchitz in \(x\), and refer to such maps as \emph{admissible}. Consequently, Carath\'eodory solutions of \eqref{eq: firstSys} exist \cite[Chapter 1]{filippov2013differential}. We tacitly assume that the solution is unique and that it exists for all time, denoting it by \(\lcro{\initT, +\infty}\ni t \mapsto x(t;{\filParam}) \in \R^d \) after suppressing the dependence on the initial state and the initial time for the sake of brevity. Without loss of generality we suppose that \(0\in\R^d\) is an equilibrium point of \eqref{eq: firstSys}; i.e., \(f(s, 0; \filParam) = 0 \) for all \(s, \filParam\). 

\begin{definition}
	\label{def:sta}
	The equilibrium point \(0\) of \eqref{eq: firstSys} is:
	\begin{enumerate}[label=\textup{(D-\alph*)}, align=left, widest=2, leftmargin=*]
		\item \label{def: UnifStab} \emph{uniformly stable} if for every \(\eps>0 \) there exists a pair \((\del, \filParam(\cdot))\), where \(\del > 0 \) and \(\filParam(\cdot)\) is an admissible map, both independent of \(t_0\) and dependent only on \(\eps\), such that 
		\[
		\normInit \leq \del \implies \norm{x(t;\filParam)} \leq \eps 	\quad \text{for all} \ t \geq t_0;
		\]
		\item \label{def: AsymStab} \textit{uniformly asymptotically stable} if for every \(\eps > 0\) there exists a pair \((\del, \ell(\cdot) )\), where \(\del> 0 \) and \(\ell(\cdot) \) is an admissible map, both independent of \(\initT \) and dependent only on \(\eps \) such that \(\normInit \leq \del \implies \norm{x(t;\filParam)} \leq \eps\) for all \(t \geq t_0\), and moreover, for each \(\eta >0 \) there exists \(T >0 \) independent of \(\initT \), such that with same admissible map \(\ell(\cdot) \) we have 
		\[
		\normInit \leq \del \implies \norm{x(t;\ell)} \leq \eta \quad \text{for all} \ t \geq t_0 + T.  
		\]
	\end{enumerate}
\end{definition}
%\begin{definition}
%	\label{def:stab}
%	The equilibrium point \(0\) of \eqref{eq: firstSys} is:
%	\begin{enumerate}[label=\textup{(\ref{def:stab}-\alph*)}, align=left, widest=2, leftmargin=*]
%		\item \label{def: UnifStab} \emph{uniformly stable} if for every \(\eps>0 \) there exists a pair \((\delta, \filParam(\cdot))\), where \(\delta > 0 \) and \(\filParam(\cdot)\) is an admissible map, both independent of \(t_0\) and dependent only on \(\eps\), such that 
%		\[
%			\normInit \leq \delta \implies \norm{x(t;\filParam)} \leq \eps 	\quad \text{for all} \ t \geq t_0;
%		\]
%	\item \label{def: AsymStable} \textit{uniformly asymptotically stable} if there exists \(r > 0\) such that for every \(\eta > 0\) there exists a pair \((T, \filParam(\cdot))\), where \(T > 0 \) and \(\filParam(\cdot)\) is an admissible map, independent of \(\initT \) and dependent only on \(\eta\), satisfying 
%		\[
%			\normInit \leq r \implies \norm{x(t;{\filParam})} \leq \eta \quad\text{for all } t \geq t_0 + T,
%		\]
%		and moreover, for each \(\eps > 0\) there exists \(\delta > 0\) independent of \(t_0\), such that with the same admissible map \(\ell(\cdot)\) we have
%		\[
%			\normInit \leq \delta \implies \norm{x(t; \filParam)} \le \eps \quad\text{for all }t\ge t_0.
%		\]
%	\end{enumerate}
%\end{definition}

\begin{remark}
	{\rm 
		{We employ slightly different versions of parametric uniform stability and uniform asymptotic stability compared to the ones predominantly found in the literature \cite[Chapter 5]{vidyasagar2002nonlinear}; this is because we must consider the dependence of the solution trajectories on the admissible map \(\filParam(\cdot)\), which plays the role of a parameter residing in some (possibly infinite-dimensional) space. A notion of parameter dependent exponential stability was introduced in \cite[Definition 1]{panteley2001uniform}, but that definition requires uniformity in the parameter in addition to the initial time, whereas we do \emph{not}. To wit, we employ the word ``uniform'' to denote uniformity in the time argument only.}
	}
\end{remark} 

%%==============================================================================
\subsection{Premise}
\label{ssec: Premise}
%%==============================================================================
For positive integers \(\stateDim\) and \(\contDim\), and an initial time \(\initT \ge 0\), consider the control affine system
\begin{equation}
\label{eq: PremiseSys}
\left\{
\begin{aligned}
& \dot{\stateTA}(t) = \drift[\stateTA(t)] + \sumFunc{k}{\contDim}\channel{k}{\stateTA(t)} u_{k}(t),\\
& \stateTA(\initT) = \bar{x}, \quad t \geq \initT,
\end{aligned}
\right.
\end{equation}
where 
\begin{enumerate}[leftmargin=*, label = \textup{(P-\roman*)}, align=left, widest=iii]
	\item \label{property:smoothness} \(\R^\stateDim \ni \xi \mapsto f(\xi) \in \R^\stateDim\) and \(\R^\stateDim \ni \xi \mapsto \channel{k}{\xi} \in \R^\stateDim\) for \(k \in \setInt[\contDim]\) are continuously differentiable maps. 
\end{enumerate}

Let \( \R^\stateDim \ni \xi \mapsto \control[](\xi) \Let \bigl(\control[k](\xi)\bigr)_{k\in\setInt[\contDim]} \in \R^m \) be a \(\contDim\)-dimensional feedback. Such feedbacks may arise from, depending on the control objective and the specific systems under consideration, the Artstein-Sontag universal formula, Jurdjevic-Quinn controllers, Lyapunov and/or dynamic programming based synthesis techniques, etc.

The application of a stabilizing feedback \(\control[]\) to \eqref{eq: PremiseSys} produces the closed-loop control system
\begin{equation}
\label{eq: normalSys}
\left\{
\begin{aligned}
& \dot{x}(t) = \drift[x(t)] + \sumFunc{k}{\contDim} \channel{k}{x(t)} \control[k](x(t)),\\
& x(\initT) = \bar{x}, \quad t \geq \initT,
\end{aligned}
\right.
\end{equation}
where we stipulate that:  
\begin{enumerate}[leftmargin=*, label = \textup{(P-\roman*)}, align=left, widest=iii]
	\setcounter{enumi}{1}
	\item \label{property:eqm} the closed-loop dynamics \eqref{eq: normalSys} has a hyperbolic,\footnote{Recall that an equilibrium point of a nonlinear system with a continuously differentiable vector field is \emph{hyperbolic} if the first order term in the Taylor's expansion of the vector field at that equilibrium point has no purely imaginary eigenvalue.} and locally asymptotically stable equilibrium at \(0\).
	\item \label{property:smoothness closed-loop} the map \(\R^\stateDim \ni \xi \mapsto \drift[\xi] + \sumFunc{k}{\contDim}\channel{k}{\xi} \control[k](\xi) \in \R^\stateDim \) is {twice} continuously differentiable { in a neighborhood of the equilibrium}, and
\end{enumerate}
As a consequence of \ref{property:eqm}, \(\drift[0] + \sumFunc{k}{\contDim}\channel{k}{0} \control[k](0) = 0\). We further assume that
\begin{enumerate}[leftmargin=*, label = \textup{(P-\roman*)}, align=left, widest=iii]
	\setcounter{enumi}{3}
	\item \label{property:simpl} \(\drift[0] = 0\), and \(\control[k](0) = 0\) for each \(k\in\setInt[\contDim]\).
\end{enumerate}
While \ref{property:simpl} is not crucial for our proofs to hold (indeed, it is sufficient to require that \(\drift[0] + \contDim\channel{k}{0} \control[k](0) = 0\) for each \(k\)), the property \ref{property:simpl} holds for a large class of standard stabilizing feedbacks including the Artstein-Sontag universal formula and the Jurdjevic-Quinn controllers to name a few, and this assumption will simplify our presentation.

As mentioned in the Introduction, the selection of the particular feedback \(\control[]\) is \emph{not} under consideration in the article at hand; we assume that there is a sufficiently large library of feedback controllers that satisfy the properties \ref{property:smoothness}-\ref{property:simpl} for a given system, from which the control designer may pick one that suits them. Instead, here we are interested in the effect of scheduling the individual components of such a stabilizing feedback in a fast round-robin fashion. This scheduling has been viewed in \cite{JurdQuinn} as a `spatial' sparsification of the control in the sense that only one of the control components is active at any given instant of time.

%%==============================================================================
\subsection{Main Results}
\label{sec:MRS}
%%==============================================================================
In this subsection we present the main results of this article. First, in \S\ref{ssec: Stability} we present the case where a family of locally asymptotically stabilizing controller for \eqref{eq: PremiseSys} is employed in a periodic fashion with constant switching frequency. Second, in \S\ref{ssec: AS} we discuss the effect of monotonically increasing the switching frequency over equispaced intervals of time. We provide a certificate of uniform stability in the sense of Definition \ref{def: UnifStab} for the closed-loop equilibrium in the former case, while in the latter case we provide a guarantee of uniform asymptotic stability in the sense of Definition \ref{def: AsymStab}. {Third, in \S\ref{ssec:BOA} we provide precise estimates of the basin of attraction of the equilibrium \(0\) under fast round-robin switching; the precise nature of these estimates will become clear from the subsequent theorems.}

%%------------------------------------------------------------------------------
\subsubsection{Uniform stability under constant switching frequency}
\label{ssec: Stability}
%%------------------------------------------------------------------------------
Fix \( \switchTime > 0 \) and \(\initT \ge 0\). Let \( \lcro{\initT, +\infty}\) be partitioned into  a countable family of disjoint half-open intervals of length \(\contDim \switchTime\) as \(\lcro{\initT, +\infty} = \bigcup_{n\in\N} \lcro{\initT+(n-1)\contDim\switchTime, \initT+n\contDim\switchTime}  \). We stipulate a periodic scheduling of the feedback controllers defined in \S\ref{ssec: Premise} over \(\lcro{\initT, \initT+\contDim\switchTime } \) with the objective of inducing temporal sparsity in the family of the \(\contDim\) controllers, and repeat the same scheme over successive such intervals. In other words, for each \(n\), we partition the interval \(\lcro{\initT+n\contDim\switchTime, \initT+(n+1)\contDim\switchTime}\) into \( \contDim\) contiguous intervals of length \(\switchTime\) each, and activate just the \(k^{th}\) controller during the \(k^{th} \) such temporal sub-interval. Accordingly, we define the following piecewise constant switching scheme: 
\begin{align}\label{eq: SwitchFunc}
	\lcro{0,+\infty} \ni t \mapsto \switchFunc(t,\switchTime) \Let 1 + \left \lfloor {\tfrac{t-\contDim\switchTime\left\lfloor{\frac{t}{\contDim\switchTime}}\right\rfloor}{\switchTime}} \right \rfloor \in \setInt[\contDim]
\end{align}
that selects, at each \(t\geq0\), the index \(\switchFunc(t, \switchTime)\) in the set \(\setInt[\contDim]\); we call \(\switchFunc\) the \textit{switching scheme} with \emph{switching time} \(\switchTime\).

Let us consider the following closed-loop system under the switching scheme \(\switchFunc\) with the feedback \(\control[](\cdot) = \bigl(\control[k](\cdot)\bigr)_{k\in\setInt[\contDim]}\) replaced by its scaled version \(\contDim\control[](\cdot) = \bigl(\contDim\control[k](\cdot)\bigr)_{k\in\setInt[\contDim]}\) in: 
\begin{equation}
\label{eq: sparseSys}
\left\{
\begin{aligned}
& \dot{\stateSS}(t) = \drift[\stateSS(t)] + \contDim \channel{\switchFunc(t-\initT,\switchTime)}{\stateSS(t)} \control[\switchFunc(t-\initT,\switchTime)](\stateSS(t)),\\
& \stateSS(\initT) = \bar{x},\quad t\ge \initT.
\end{aligned} 
\right.
\end{equation}
In the light of \ref{property:simpl}, it is readily observed that the point \(0\in\R^d\) is an equilibrium point of the switched system \eqref{eq: sparseSys}.

\begin{remark}\label{rem: Amplify}
	{\rm We have scaled each of the feedback functions \(\control[k](\cdot)\) by a factor of \(\contDim\) in \eqref{eq: sparseSys}. Since the switching scheme \(\switchFunc\) permits only one \(\control[k]\) to be active `on' at a given time, this extra factor of \(\contDim\), roughly speaking, extracts sufficient ``energy'' from the system during its `on' stage to compensate for the inactivity during its `off' period in order to stabilize the equilibrium point.}
\end{remark}

Theorem \ref{thm: MainResult} below, {whose proof is deferred to Appendix \ref{app_lab: Main Result}, asserts that under the round-robin periodic scheduling with constant switching frequency, the equilibrium point \(0\in\R^d\) of the closed-loop system \eqref{eq: sparseSys} is uniformly stable in the sense of Definition \ref{def: UnifStab}.
	
	\begin{theorem}
		\label{thm: MainResult}
		Consider the control system \eqref{eq: PremiseSys} along with its associated data \ref{property:smoothness}. Suppose that we pick a feedback \(\control[](\cdot)\) such that the closed-loop system \eqref{eq: normalSys} satisfies \ref{property:eqm} -\ref{property:simpl}. Then the equilibrium \(0 \in \R^d \) of \eqref{eq: sparseSys} is uniformly stable in the sense of Definition \ref{def: UnifStab} ({considering the switching time as a parameter}). %, where the constant \(\switchTime \) is playing the role of the time-dependent map \(\filParam(\cdot) \).
		Moreover, for every \(\eps> 0 \), there exist \(\del, \switchTimeS > 0 \), both independent of \(\initT \), such that for all \(\switchTime \in \lorc{0, \switchTimeS} \) and with the switching scheme \(\switchFunc(\cdot,\tau) \) defined in \eqref{eq: SwitchFunc}, the solution of \eqref{eq: sparseSys} satisfies 
		\[
		\normInit \leq \del \implies \norm{y(t)} \leq \eps\quad\text{for all } t \geq \initT.
		\] 
	\end{theorem}
	
	This theorem caters to the `fast switching' regime in the switched systems literature, an area where there have been far fewer investigations compared to the more well-known `slow switching' regime treated extensively in the standard textbook \cite[Chapters 3, 4]{LibBook} and more recently surveyed in \cite{ref:KunCha-17}. Intuition suggests that Theorem \ref{thm: MainResult} pertains to a type of averaging on a fast time-scale. Indeed, the chattering lemma, popularly employed for various types of convexification in the context of optimal control theory, is the key tool needed in our proofs. This state of affairs points to the possibility of employing, more generally, Young measures in the context of averaging along the lines of \cite{ref:Art-08}; this is currently under investigation and will be reported in subsequent articles.
	
	%\begin{remark}
	%	Note that the switching time \(\switchTimeS \) for which the bound in \eqref{eq: unifBoundMR} holds is, in fact, very conservative. This is because it is derived assuming worst case behavior of the system trajectories which correspond to the case when the vector fields assume maximal possible values (depending on the Lipchitz constant). In general, depending on the system dynamics one may get uniform stability even at higher switching time (or correspondingly lower switching frequency). 
	%\end{remark}
	%\begin{remark}
	%	The choice of \(\pos[1] \) is free in the interval \(\loro{0,1} \). This choice of \(\pos[1] \) governs the upper bound on \(\gamBar \). The higher the \(\pos[1] \), the upper bound on \(\gamBar \) reduces. Further, the choice of the switching time \(\switchTime \) is dependent upon \(\gamBar\). If we control engineer wants that the trajectories of \eqref{eq: sparseSys} should follow the trajectories of \eqref{eq: normalSys} closely then smaller switching time (or higher switching frequency) should be employed. 
	%\end{remark}
	
	%%------------------------------------------------------------------------------
	\subsubsection{Uniform asymptotic stability under increasing switching frequencies}
	\label{ssec: AS}
	%%------------------------------------------------------------------------------
	Theorem \ref{thm: MainResult} dealt with uniform stability of the round-robin periodically scheduled (switched) system \eqref{eq: sparseSys} under a constant switching frequency. However, in control engineering we are often interested in ensuring asymptotic convergence of system trajectories to the equilibrium point.
	% Moreover, we have assumed that the trajectories of \eqref{eq: normalSys} are locally asymptotically stable to the equilibrium and therefore it is obvious to look for sufficient conditions for asymptotic stability of the trajectories of corresponding sparse system \eqref{eq: sparseSys}.
	It turns out (as we shall assert in Theorem \ref{thm: AS_stable}), that increasing the switching frequency periodically so that the frequency becomes unbounded asymptotically in time guarantees uniform asymptotic stability of the equilibrium point of the resulting switched system.
	
	Fix \(\initT \ge 0 \) and \(\timeDown[] > 0 \). Let \(\lcro{0, +\infty} \ni t \mapsto \switchTimeAS(t) \in \loro{0, +\infty}\) be a piecewise constant non-increasing function satisfying \(\switchTimeAS(t) \xrightarrow[t\to+\infty]{} 0\). If \(\tau':\lcro{0, +\infty}\lra\R\) is another piecewise constant non-increasing function satisfying \(\tau'(t)\xrightarrow[t\to+\infty]{} 0\), we say that \(\switchTimeAS\) \emph{dominates} \(\tau'\) if \(\switchTimeAS(t) \le \tau'(t)\) for each \(t\).\footnote{The word \emph{domination} here sounds more appropriate when we note that the `frequency' of \(\switchTimeAS\) compared to \(\tau'\) is indeed higher.} We define a new (time-varying) \(\switchTimeAS\)-dependent switching scheme
	\begin{align}\label{eq: switchRuleAS}
		\lcro{0, +\infty} \ni t \mapsto \switchFuncAs(t; \switchTimeAS) \Let 1 + \Bigl\lfloor {\tfrac{t-\contDim\switchTimeAS(t)\left\lfloor{\frac{t}{\contDim\switchTimeAS(t)}}\right\rfloor}{\switchTimeAS(t)}} \Bigr\rfloor \in \setInt[\contDim].
	\end{align}
	The interval \(\lcro{\initT, +\infty} \) is partitioned, depending on \(\switchTimeAS(\cdot)\), into a countable family of disjoint half-open intervals of length \(\timeDown[]\) as \(\lcro{\initT, +\infty} = \bigcup_{n\in\N} J_n\), where \( J_n \Let  \lcro{\initT + (n-1)\timeDown[], \initT + n\timeDown[]} \); by definition, on each interval \(J_n\) the function \(\switchTimeAS(\cdot) \) is held constant. 
	
	Here we are interested in establishing the asymptotic convergence of the following closed-loop system derived from \eqref{eq: normalSys} under the switching scheme \(\switchFuncAs\) with the feedback \(\control[](\cdot) = \bigl(\control[k](\cdot)\bigr)_{k\in\setInt[\contDim]}\) replaced by its scaled version \(\contDim\control[](\cdot) = \bigl(\contDim\control[k](\cdot)\bigr)_{k\in\setInt[\contDim]}\): 
	\begin{equation}
	\label{eq: sparseSysAS}
	\left\{
	\begin{aligned}
	& \dot{\stateAS}(t) = \drift[\stateAS(t)] + \contDim \,\channel{\switchFuncAs(t-\initT;\switchTimeAS)}{\stateAS(t)} \,\control[\switchFuncAs(t-\initT;\switchTimeAS)](\stateAS(t)),\\
	& \stateAS(\initT) = \stinit,\quad t\ge \initT,
	\end{aligned} 
	\right.
	\end{equation}
	where \(\stinit \) is the same initial vector as in \eqref{eq: normalSys}. 
	
	Against the preceding backdrop, we present the following theorem whose proof is deferred to Appendix \ref{app_label: ASstable}. It guarantees the existence of a piecewise constant non-increasing function \(\switchTimeAS \) that encodes the switching time information, such that the equilibrium for the closed-loop system \eqref{eq: sparseSysAS} is locally uniformly asymptotically stable  in the sense of Definition \ref{def: AsymStab}.  
	
	\begin{theorem}
		\label{thm: AS_stable}
		Consider the control system \eqref{eq: PremiseSys} along with its associated data \ref{property:smoothness}. Suppose that we pick a feedback \(\control[](\cdot)\) such that the closed-loop system \eqref{eq: normalSys} satisfies \ref{property:eqm}-\ref{property:simpl}. Then the equilibrium of \(0 \in \R^d \) of \eqref{eq: sparseSysAS} is locally uniformly asymptotically stable in the sense of Definition \ref{def: AsymStab} ({considering the switching time as a parameter}). % where the map \(\switchTimeAS \) is playing the role of the time-dependent map \(\filParam(\cdot) \).
		Moreover, for every \(\eps > 0 \), there exists a pair \((\del, \switchTime'(\cdot)) \), independent of \(\initT \), where \(\del > 0 \) and \(\lcro{\initT, +\infty} \ni t \mapsto \switchTime'(t) \in \loro{0, +\infty} \) is a piecewise constant non-increasing map, such that
		\begin{itemize}[label=\(\circ\)]
			\item for every \(\switchTimeAS\) that dominates \(\tau'\), with the switching scheme \(\switchFuncAs(\cdot;\switchTimeAS) \) defined in \eqref{eq: switchRuleAS}, we have 
			\[
			\normInit \leq \del \implies \norm{z(t)} \leq \eps \quad \text{for all}\ t \geq \initT, 
			\]
			and
			\item for every \(\eta > 0 \), there exists \(\tilde{T} > 0 \), independent of \(\initT \), such that with the same switching scheme \(\switchFuncAs(\cdot;\switchTimeAS)\), we have 
			\[
			\normInit \leq \del \implies \norm{z(t)} \leq \eta \quad \text{for all}\ t \geq \initT+\tilde{T}. 
			\]
		\end{itemize}
	\end{theorem}
	
	\begin{remark}
		\label{rem: allChannels}
		{\rm The underlying premise of the preceding discussion is that at any instant of time only \emph{one} out of \(\contDim \) feedbacks is permitted to stay `on'. However, our analysis carries over to the situation when there is an ensemble of control-affine nonlinear systems, and each control input \(u_i\) is itself a multi-variable control catering to its own stand-alone system; we construct the aggregate high-dimensional control system by stacking the individual systems in a natural way to obtain the aggregated high-dimensional drift and control vector fields.}%Although there is no a priori coupling between the individual systems at the stage of this definition, a nominal (feedback) controller for this aggregated system may well introduce such coupling terms via the resulting feedbacks. An identical analysis as ours may be carried hereafter by scheduling the family of inputs \((u_i)\) as discussed above.
	\end{remark}

	\begin{remark}\label{rem: JurdQuinn}
		{\rm Let us highlight a few differences between \cite{JurdQuinn} and our results. \cite[Theorem 1.1]{JurdQuinn} does \emph{not} apply to any linear multivariable control system having an unstable system matrix, whereas Theorems \ref{thm: MainResult} and \ref{thm: AS_stable} do, and we illustrate the application of the latter two results to a linear multivariable control system in \S\ref{sec: LinModel}. For nonlinear control systems, the situation is more complicated: \cite[Theorem 1.1]{JurdQuinn} does not stipulate that the Lyapunov function \(V\) in that article is locally positive definite. It is possible, therefore, for \(V\) to vanish on a compactum having non-empty interior. However, both Theorem 1.1 and Corollary 1.2 of the aforementioned article expressly require that either an invariant subset (denoted
			by $\Omega$ there) of $\mathcal{Z}$ in op.\ cit.\ or the origin be attractive under zero control, in which case the theorem asserts the convergence of trajectories under sparsification to this set. When the set \(\mathcal Z\) in op.\ cit.\ is specialized to \(\{0\}\) to conform to our setting, the requirement of \(\mathcal Z = \{0\}\) being attractive continues to be an integral part of \cite[Corollary 1.2]{JurdQuinn}. In contrast, the point \(0 \in \R^d\) is not attractive for the `zero-input' system in our setting.}
	\end{remark}
	
	\begin{remark}\label{rem: BacMaz}
		{\rm Let us highlight a few differences between \cite{ref:BacMaz-10} and our results. The premise of \cite{ref:BacMaz-10} is a finite family \((f_i)_{i\in\mathcal N}\) of smooth vector fields on \(\R^d\) that vanish at \(0\in\R^d\), and such that for some convex combination \((\alpha_i)_{i\in\mathcal N}\) of the family \((f_i)_{i\in\mathcal N}\), \(0\in\R^d\) is an asymptotically stable equilibrium of the system \(\dot x = \sum_{i\in\mathcal N} \alpha_i f_i(x)\). For such a family, \cite[Theorem 1]{ref:BacMaz-10} constructs a time- and state-dependent multiplexing/switching rule \(\sigma\) that stabilizes the point \(0\in\R^d\) for the switched system \(\dot x(t) = f_{\sigma(t)}(x(t))\). The process of designing \(\sigma\) in op.\ cit.\ is rather tedious and difficult for \(\contDim \geq 3\) because the techniques rely on the Baker-Campbell-Hausdorff formula, has no relationship with fast switching, and as pointed out in Remark 3 of op.\ cit., \(\sigma\) may be aperiodic and may depend on the distance from \(0\) of the states. One way to apply \cite[Theorem 1]{ref:BacMaz-10} in our context would be to pick \(\mathcal N = \{1, \ldots, m\}\) and \(f_i = f + m g_i\) for \(i = 1, \ldots, m\); then (the proof of) this theorem would provide \emph{a particular} time- and space-dependent \(\sigma\) such that \(\dot x(t) = f(x(t)) + m g_{\sigma(t)}(x(t))\) is asymptotically stable. In contrast, our results show that sufficiently fast and simple round-robin sparsification of \emph{any} stabilizing multivariable feedback controller preserves stability; our scheme relies on fast switching between the selected family of scalar feedbacks, and our results also describe the behaviour of the controlled system under increasingly faster switching. Our results and the ones in \cite{ref:BacMaz-10} are, therefore, complementary.}
	\end{remark}
	
	\begin{remark}\label{rem: NCSLit}
		{\rm The articles \cite{Nesic2004Automatica,Nesic2004TAC,walsh2001asymptotic} in the area of networked control systems (NCS) and the more recent work \cite{wang2012input} are relevant to our problem, but the hypotheses employed in these articles are {\emph{stricter} than ours when restricted to our setting of control affine systems; consequently, our results are \emph{novel}}. In particular, all these cited works employ different types of quantitative (and sometimes restrictive) estimates of the form \cite[Definitions 4 and 5]{wang2012input}, which permit them to employ high-level tools. We do \emph{not} employ any quantitative estimates and rely solely on a delicate analysis of continuity of solutions to Carath\'eodory equations derived from the techniques exposed in \cite[Chapter 1]{filippov2013differential}. This is the subject of section A.2. Notice that while articles along the lines of \cite{wang2012input} in the NCS literature discuss time-averaging of vector fields (for hybrid systems in some articles), we rely on a parameter dependent averaging of a different kind. Therefore in Theorem A.3 (borrowed from \cite{filippov2013differential}), the weak convergence of the vector fields assumed, in equation (A.5), is dependent on the `parameter' sequence (\( l_i \)). The articles \cite{Nesic2004TAC,Nesic2004Automatica} in the NCS literature do consider the problem of retaining some properties (such as \(\mathcal{L}_p \) stability or input to state stability) under a version of sparsification that occurs due to transmission protocols in a networked system. In that sense there is an underlying similarity. {However, the phrase ``round-robin'' has been defined differently in these articles: In their context the transmission protocol facilitates access to one set of nodes/agents whose outputs and controls are updated (typically to current true values) while the outputs/controls of other nodes/agents are \emph{retained at the old values}; their round-robin protocol ensures that this access is granted by cycling through all nodes/agents periodically. This is in stark contrast to our situation: While the controllers are still cycled periodically, at each `switch', the control of one node (if thought of as a networked system) is updated to the nominal value while the control at the other nodes are changed to \emph{zero} (as opposed to their previous nominal values). Therefore, proving that our `protocol' is uniformly globally asymptotically stable via existing techniques appears to be difficult since while the control error on one node improves, the error on all other nodes deteriorates significantly.}}
	\end{remark}
	
	\subsubsection{Basins of attraction}
	\label{ssec:BOA}
	
	{
		A natural question that arises from the preceding discussion is: How does the basin of attraction of the asymptotically stable equilibrium point \(0\in\R^d\) of the nominal system change under the proposed round-robin scheduling? We shall provide answers to this question below. Let us denote the basin of attraction of the equilibrium in \eqref{eq: normalSys} by \(\roa{cl}\), and in order to transparently carry the dependence of the initial conditions and the switching scheme, we shall denote the trajectories of \eqref{eq: sparseSysAS} by \(\stateAS(\cdot;\stinit,\switchFuncAs(\cdot;\switchTimeAS),\initT) \) until the end of this section.
		
		The following theorem, whose proof is deferred to Appendix \ref{ssec: NonUniform}, asserts that the basin of attraction of the round-robin scheduled system is identical to the nominal system under appropriate scheduling:
		
		\begin{theorem}\label{thm: NonuniformResult}
			Consider the control system \eqref{eq: PremiseSys} along with its associated data \ref{property:smoothness}. Suppose that we pick a feedback \(\control[](\cdot)\) such that the closed-loop system \eqref{eq: normalSys} satisfies \ref{property:eqm}-\ref{property:simpl}. Then for every point \(\stinit \in  \roa{cl}\) there exists a switching time \(\switchTime'(\cdot;\stinit)\) (depending on the initial condition \(\bar{x}\)) such that for any \(\switchTimeAS\) that dominates \(\switchTime'\),  \(\stateAS(t;\stinit,\switchFuncAs(\cdot;\switchTimeAS), \initT)\xrightarrow[t\ra+\infty]{} 0\). 
		\end{theorem}
		
		We emphasize that the switching time in the previous theorem depends on the initial condition. One can circumvent this dependence by restricting to interior of the set \(\roa{cl}\) in the following way. First we recall some machinery: Since the equilibrium \(0\in \R^d\) of \eqref{eq: normalSys} is asymptotically stable, the converse Lyapunov theorem \cite[Chapter 5]{vidyasagar2002nonlinear} guarantees that there exists a continuously differentiable positive definite function \(\lyap: \roa{cl}\ra \lcro{0,+\infty}\) such that \(V(x)\xrightarrow[x\ra\bdd\roa{cl}]{} +\infty\), and it decreases strictly along the trajectories of \eqref{eq: normalSys}. For any \(\radCompact>0\) define \(\compactSet{\radCompact} \defas \{x\in\roa{cl}: \lyap(x) \leq \radCompact\}\) to be the \(\radCompact\)-sublevel set of \(\lyap\). This set is compact and invariant for every \(\radCompact > 0\) \cite[Chapter 5]{vidyasagar2002nonlinear}. We have the following theorem, whose proof is deferred to Appendix \ref{ssec: Uniform}: 
		
		\begin{theorem}\label{thm: RegionOfAttractionUniformResult}
			%Consider the control system \eqref{eq: PremiseSys} along with its associated data \ref{property:smoothness}. Suppose that we pick a feedback \(\control[](\cdot)\) such that the closed-loop system \eqref{eq: normalSys} satisfies \ref{property:smoothness closed-loop}-\ref{property:simpl}. Let \(\contract \in \loro{0,1}\). 
			Consider the control system \eqref{eq: PremiseSys} along with its associated data \ref{property:smoothness}. Suppose that we pick a feedback \(\control[](\cdot)\) such that the closed-loop system \eqref{eq: normalSys} satisfies \ref{property:eqm}-\ref{property:simpl}. Then for every \(\radCompact > 0\) there exists a switching time \(\switchTime'(\cdot;\radCompact)\) (depending on \(\radCompact\)) such that for any \(\stinit \in \compactSet{\radCompact}\) and any \(\switchTimeAS\) that dominates \(\switchTime'\), we have \(\stateAS(\cdot;\stinit,\switchFuncAs(\cdot,\switchTimeAS),\initT) \xrightarrow[t\ra+\infty]{} 0\).
		\end{theorem}
		%Informally, both of the results say that under sufficiently high switching frequency   `difference' between \roa{cl}\ and \roa{rr}\ can be made arbitrarily small.  
	}
	
	%\todo[inline, color=Yellow!20]{Srikant, I don't understand a word of this remark:}
	%{\color{BrickRed} 
	%\begin{remark}
	%	{\rm The idea of ``approximating'' the time varying system with that of time invariant system is instrumental in our analysis \S\ref{app: PTheorem}. In spirit this idea is similar to that in \cite{wang2012input} where the time varying system is approximated by suitably defined time invariant systems. However, their time average of time varying system approaches that of  time invariant systems (refered as weak and strong average systems \cite[Definition 4,5]{wang2012input}) as the horizon increases, which is in contrast to our work where we define the invariant system to be the limit of time varying system in limit of a parameter.}
	%\end{remark}
	%}

%%%%%%%%%%%%%%%%%%%%%%%%%%%%%%%%%%%%%%%%%%%%%%%%%%%%%%%%%%%%%%%%%%%%%%%%%%%%%%%%
\section{Numerical Experiments}
\label{sec: NumExp}
%%%%%%%%%%%%%%%%%%%%%%%%%%%%%%%%%%%%%%%%%%%%%%%%%%%%%%%%%%%%%%%%%%%%%%%%%%%%%%%% 
In this section we present comprehensive numerical studies of two control systems, namely, a coupled inverted pendulum on a cart and a mass-spring system; and a spacecraft attitude control system when the control inputs are scheduled in a periodic fashion as discussed above. In the former example,	we employ a constant switching rate based scheduling (discussed in \S\ref{ssec: Stability}) while in the latter we use time-varying scheduling with the switching rate monotonically increasing and diverging with time (discussed in \S\ref{ssec: AS}). 
%%==============================================================================
\subsection{Linearized model of an inverted pendulum on a cart coupled with a mass-spring system}
\label{sec: LinModel}
%%==============================================================================

We start with a contrived example derived from the (Taylor's) linearization of an inverted pendulum on a cart about its upright unstable equilibrium point and a mass-spring system. This example illustrates several features of Theorem \ref{thm: MainResult}, including the fact that for a fixed stabilizing feedback, 
\begin{itemize}
	\item it is possible to ensure uniform stability under fast round-robin switching,
	\item slow switching may lead to instability, and
	\item omitting the scaling of the control actions by \(m\) may lead to instability when the corresponding scaled control actions gives uniform stability.
\end{itemize}

The control system at hand can be represented as follows: 
\begin{align}\label{eq: crossCouple}
	\dot{x}(t) = \underbrace{\begin{pmatrix}
			{0} & {1} & {0}& {0}& 0 & 0 \\
			{0} & {0} & {\frac{-mg}{M}} & {0} & 0 &0\\
			{0} & {0} & {0} & {1} & 0 & 0 \\
			{0} & {0} & {\frac{(m+M)g}{LM}} & {0} & 0&0 \\
			0 & 0 & 0 & 0 & {0} & {1} \\
			0 & 0 & 0 & 0 & {-1} & {0}
	\end{pmatrix}}_{\teL A} x(t)
	+ \underbrace{\begin{pmatrix}
			{0} & 0 \\
			{\frac{1}{M}} & {0.1}\\
			{0} & 0 \\
			{\frac{-1}{LM}} & {0.1} \\
			0 & {0} \\
			{0.1} & {1}	 
	\end{pmatrix}}_{\teL B}\begin{pmatrix}
		{u_1}(t) \\ {u_2}(t)
	\end{pmatrix},
\end{align}
where $x \Let (x_{IP}, x_{MS})\in\R^6$, with $x_{IP} \in \R^4$ denoting the states of the inverted pendulum on a cart, and $x_{MS} \in \R^2$ representing the states of the mass-spring system. In \eqref{eq: crossCouple} the quantity \(M \) denotes the mass of the cart, \(m \) denotes the mass of the pendulum, \(L \) is the length of the pendulum, and \(g\) refers to the acceleration due to gravity.

The drift vector field in \eqref{eq: crossCouple} corresponding to the inverted pendulum is unstable due to the fact that the upright vertical position of the inverted pendulum is unstable. The two systems are coupled via the matrix \(B \) in \eqref{eq: crossCouple}, and this coupling has been artificially introduced for illustrative purposes only and is not motivated by physical considerations. Due to this coupling, the control input that corresponds ideally to the mass-spring system may contribute to a detrimental behaviour of the pendulum during its `on' periods, and conversely.

Recall that if there exists a matrix $K$ such that the matrix $(A-BK)$ is Hurwitz then the feedback controller $u = -Kx$ when applied to \eqref{eq: crossCouple} makes the corresponding closed-loop system globally exponentially stable at the origin (hence it is asymptotically stable as well). We pick any such gain \( K \), and employ this feedback via the switching scheme \eqref{eq: SwitchFunc} to arrive at the analogue of \eqref{eq: sparseSys} for the present example. The data employed in our experiments is presented in Table \ref{tab: Qtt1}.
\begin{table}
	\begin{center}
		\begin{tabular}{ @{}ll@{} }
			\toprule
			Quantity & Values \\
			\midrule
			Time of simulation & 50 \si{\second} \\  
			Switching time & 0.5 \si{\second}\\
			ODE solver & \(4\)-th order RK \\
			Step length & 0.05 \si{\second}\\
			Mass of the cart & 3 \si{\kilogram}\\
			Mass of the pendulum & 0.25 \si{\kilogram}\\
			Acceleration due to gravity & 9.81 \si{\kilogram\meter\second^{-2}}\\
			Initial condition & \((0, \frac{\pi}{10}, 0, 0, 1, 1.05)^\top\) \\
			\bottomrule
		\end{tabular}
		\caption{Parameters used in \eqref{eq: crossCouple}.}
		\label{tab: Qtt1}
	\end{center}
\end{table}

As per our discussion in the preceding sections, in the periodic scheduling we keep all but one feedbacks in the `off' mode (see Fig.\ \ref{fig: Control}). The plant in the current control system has at least one unstable direction due to the unstable behavior of the pendulum about its upright vertical state; this leads to the divergence of the trajectories during the off-phase of the control actuation for the corresponding states. Some of the system trajectories in Fig.\ \ref{fig: InvPend} and Fig.\ \ref{fig: SHO} exhibit fluctuations precisely due to this feature.

Apart from the fast round-robin scheduling, the role of the amplification of the control action (discussed in Remark \ref{rem: Amplify}) proves to be important for uniform stability. Indeed, if we remove this amplification factor, then the system trajectories may, in fact, diverge, as illustrated in Fig.\ \ref{fig: blowingtraj}. 

Theorem \ref{thm: MainResult} asserts that the equilibrium point of the system under fast periodic scheduling of its controllers is uniformly stable. However, it is interesting to note that the trajectories in this particular example appear to be exponentially stable (see Fig.\ \ref{fig: norm_x}). A formal proof of the same is not available at this juncture.  

\begin{figure}
	\centering
	\includegraphics[width = 0.6\textwidth]{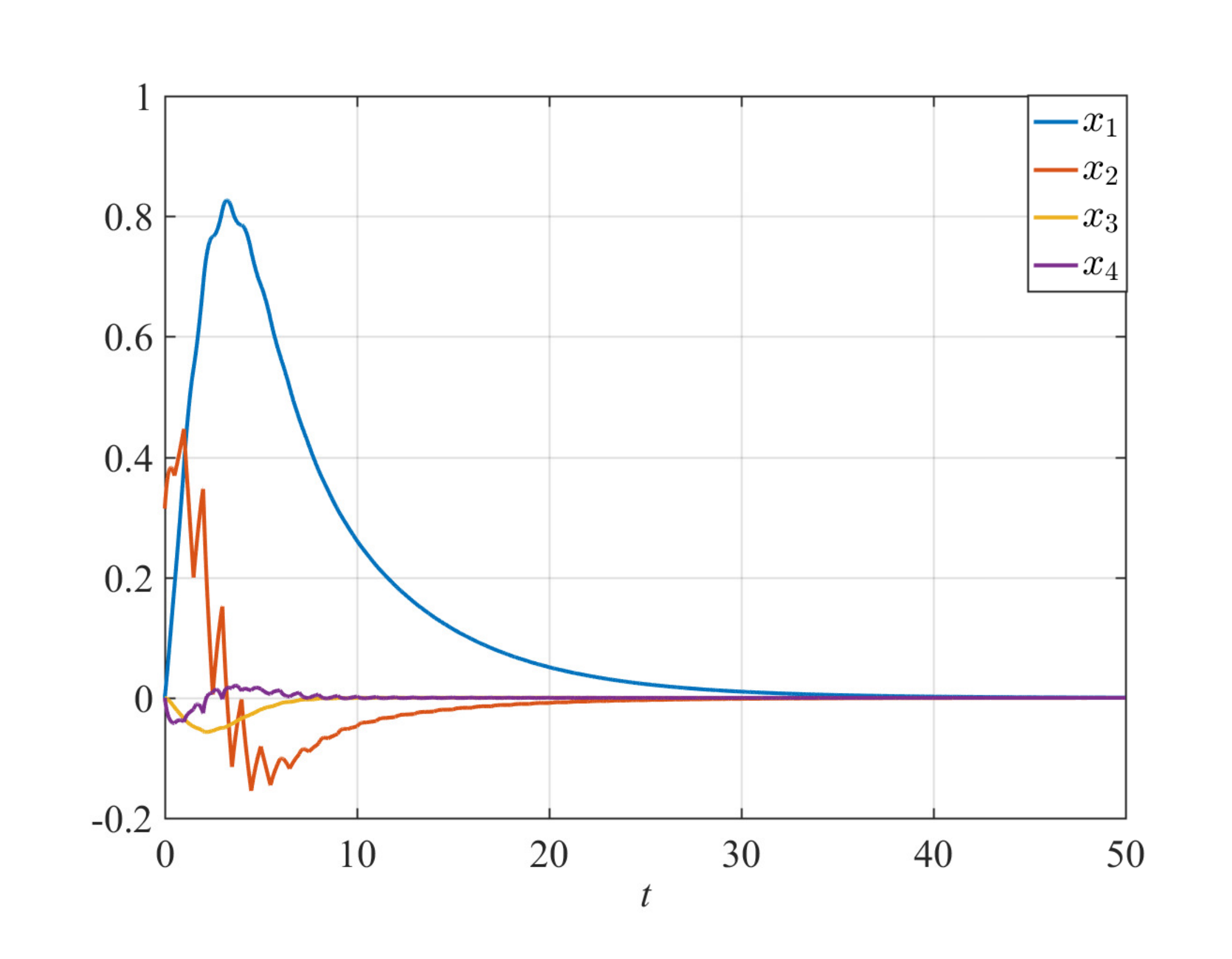}
	\caption{The time evolution of the states of the inverted pendulum on a cart.}
	\label{fig: InvPend}
\end{figure}
\begin{figure}
	\centering
	\includegraphics[width = 0.6\textwidth,keepaspectratio]{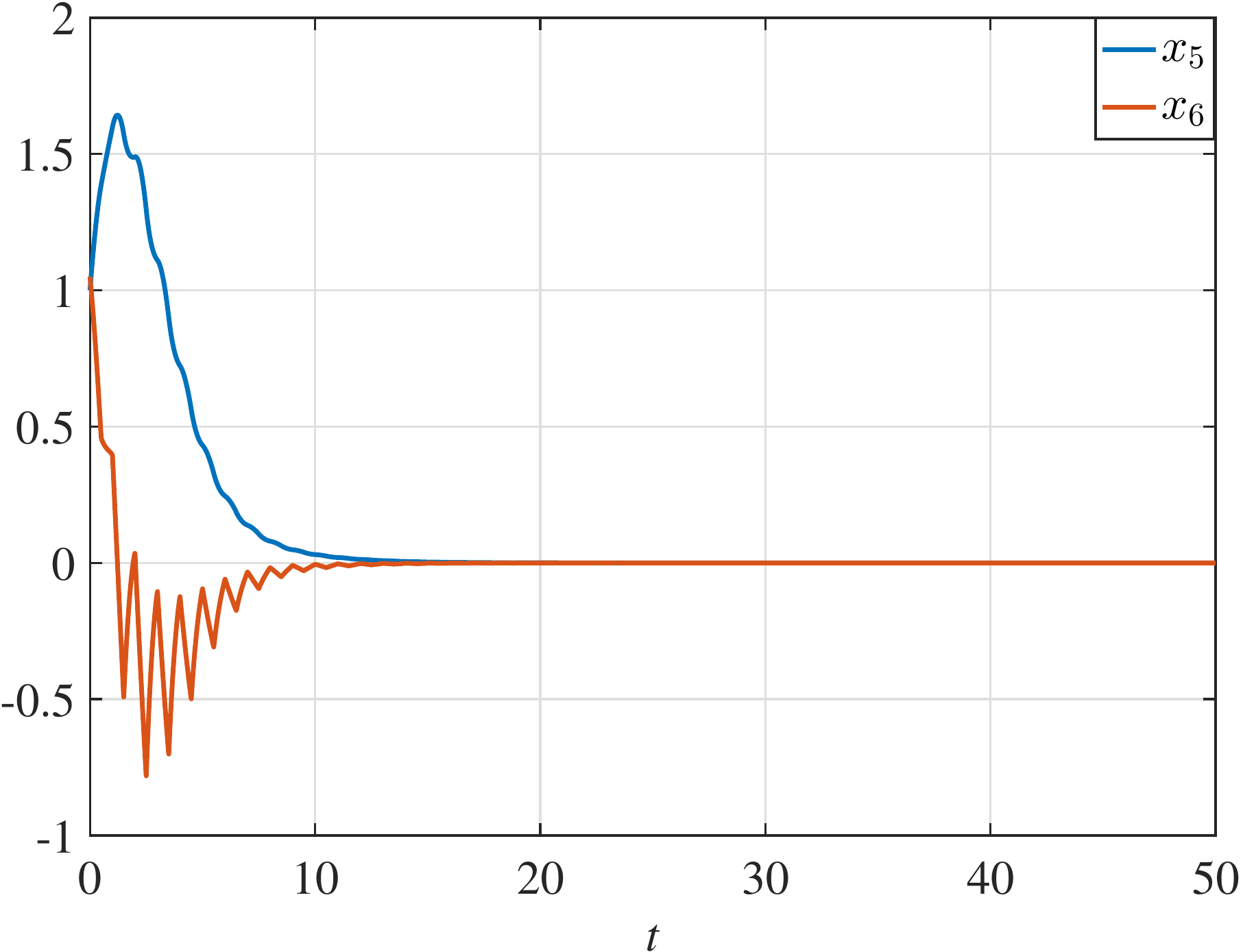}
	\caption{The time evolution of the states of the mass-spring system.}
	\label{fig: SHO}
\end{figure}
\begin{figure}[h!]
	\centering
	\includegraphics[width = 0.6\textwidth,keepaspectratio]{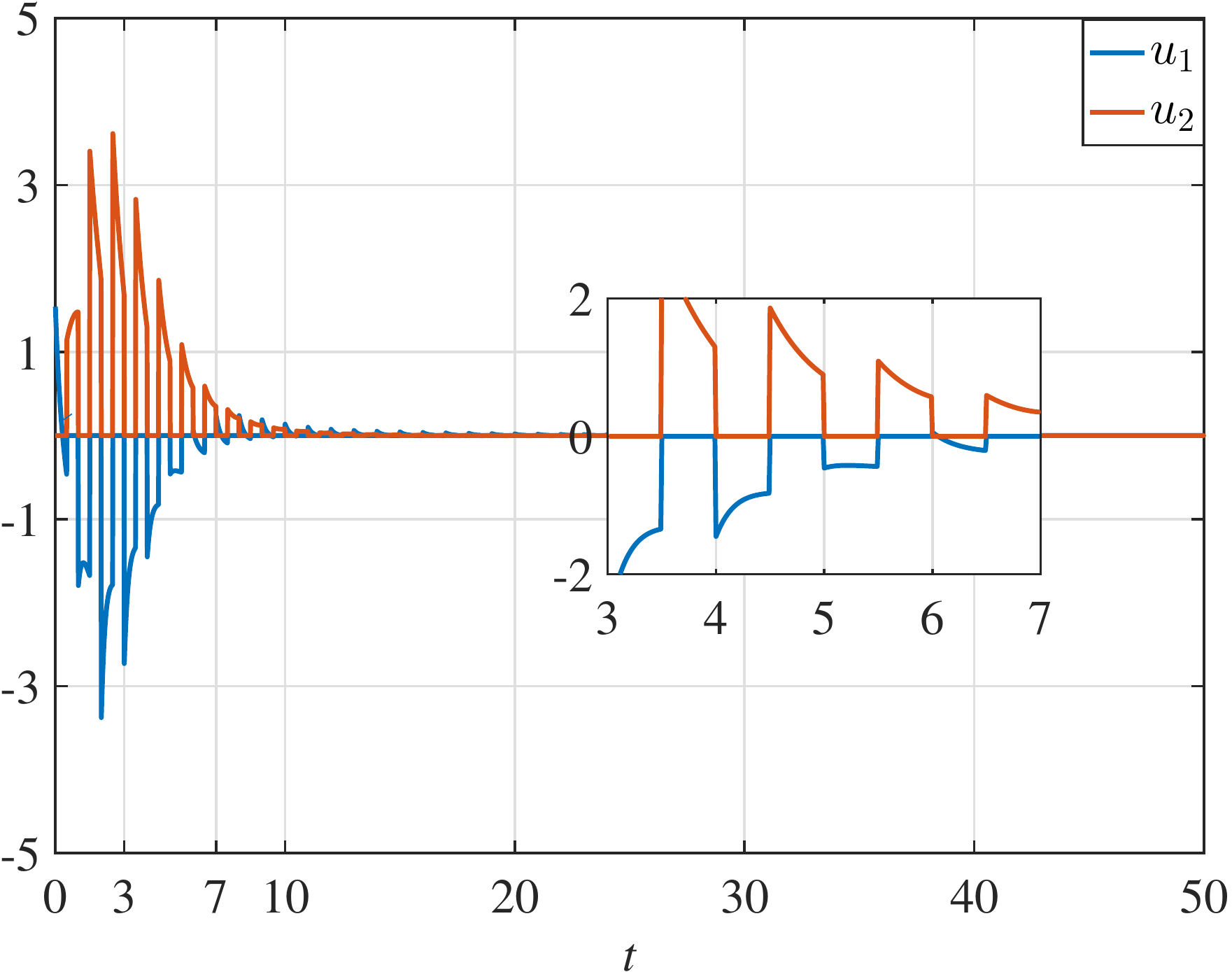}
	\caption{Control trajectories for the inverted pendulum (in blue) and the mass-spring system (in red)}
	\label{fig: Control}
\end{figure}
\begin{figure}
	\centering
	\includegraphics[width = 0.6\textwidth,keepaspectratio]{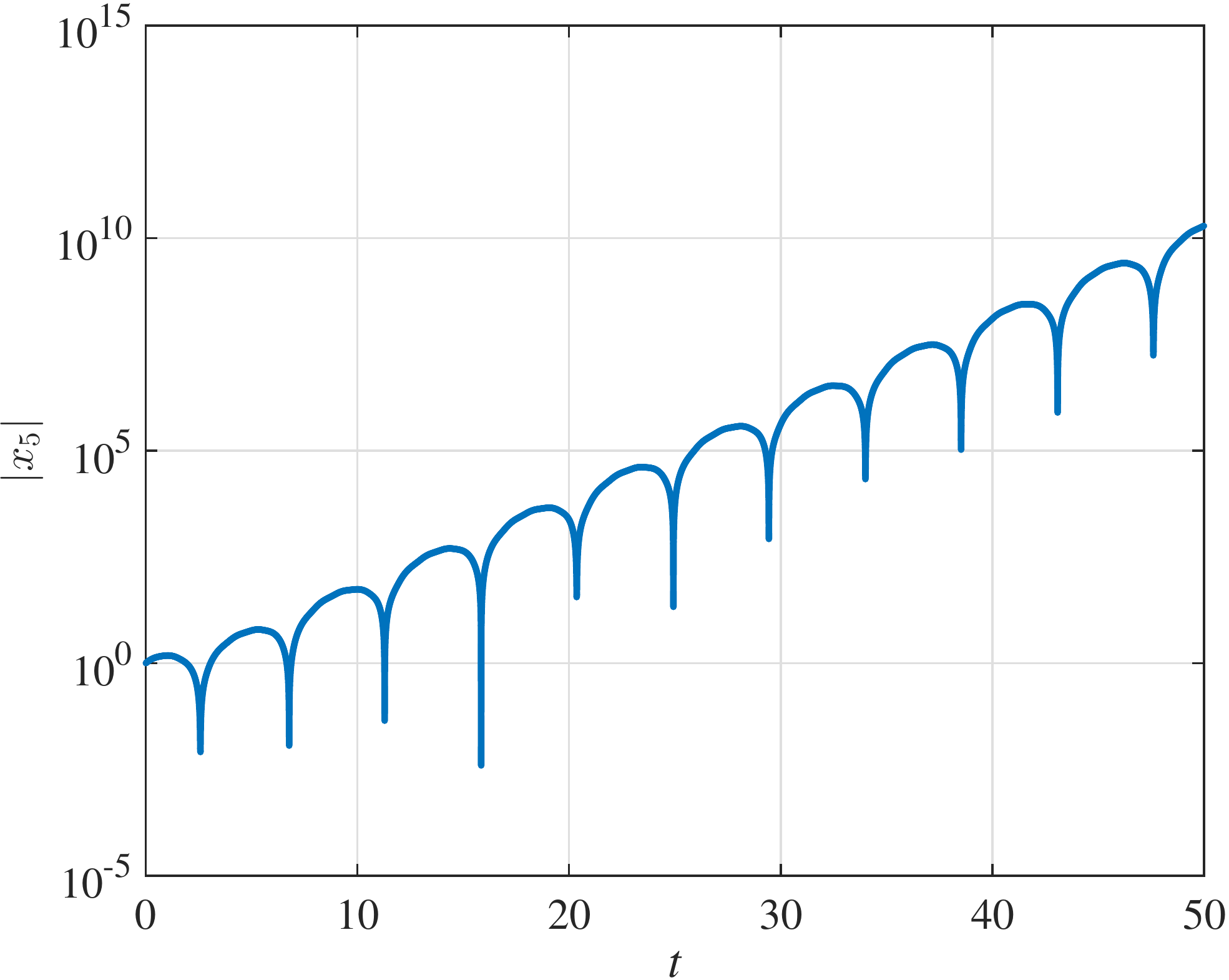}
	\caption{Time evolution of \(|x_5(\cdot)| \) in \eqref{eq: crossCouple} when the multiplicative factor of \(m\) is omitted; exponential divergence is evident.}
	\label{fig: blowingtraj}
\end{figure}     
\begin{figure}
	\centering
	\includegraphics[width = 0.6\textwidth,keepaspectratio]{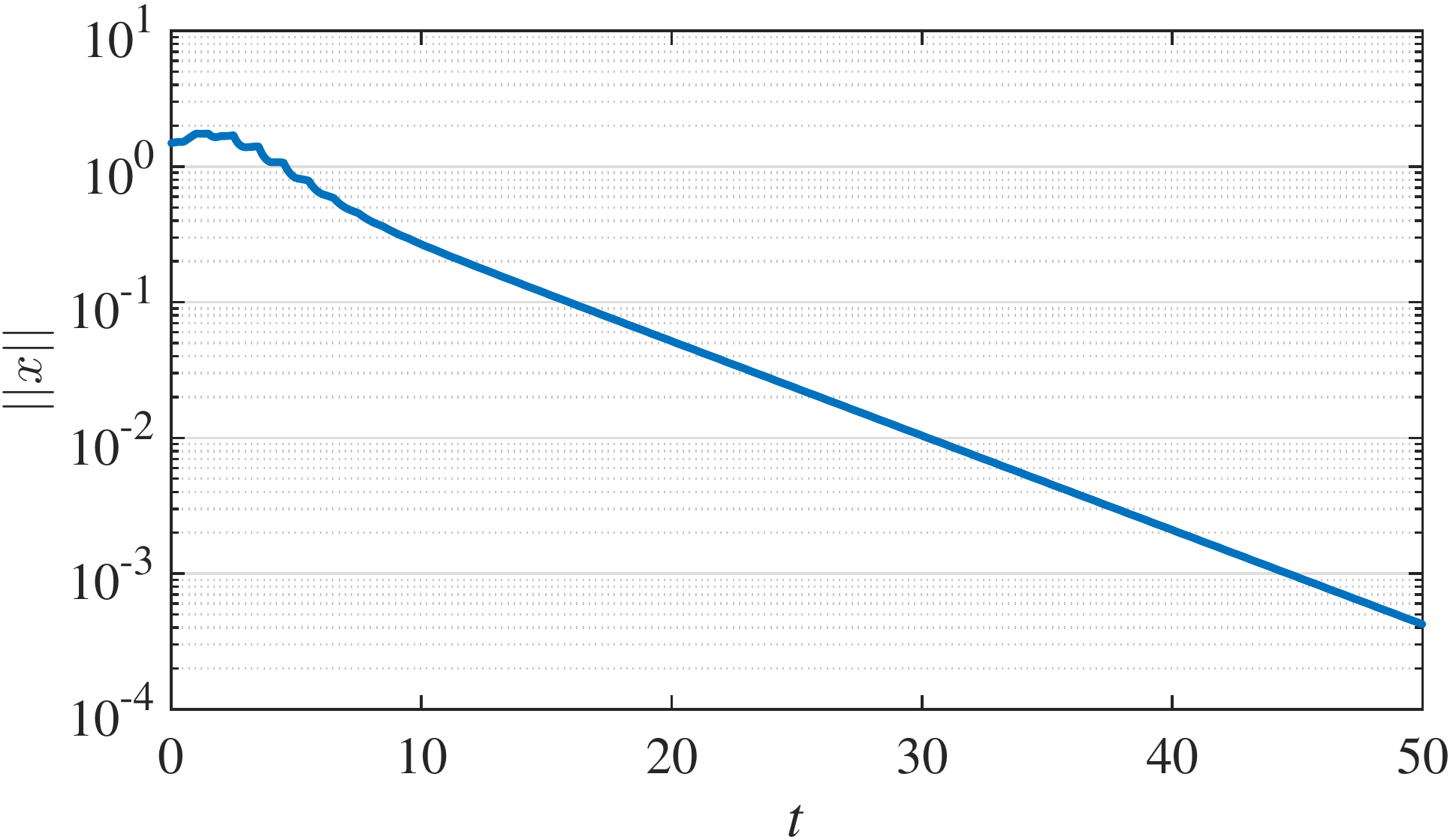}
	\caption{Time evolution of $\|x(\cdot)\|$ on a $\log$ scale for \eqref{eq: crossCouple}; exponential decay is evident.}
	\label{fig: norm_x}
\end{figure}   

%%==============================================================================
\subsection{Attitude control of spacecraft}
\label{sec: Spacecraft}
%%==============================================================================
Spacecraft are typically laden with instruments such as antennas, solar panels, etc., which must be pointed in specific directions in space for optimal performance of these instruments. The attitude (orientation) control of spacecraft is, therefore, of paramount importance today; more details on this control problem may be found in \cite{ref:wen-99}. In this section we deploy an asymptotically stabilizing controller for spacecraft attitude in the round-robin fashion with time-varying switching rates. This example illustrates several features of our main results, including
\begin{itemize}
	\item asymptotic stability under monotonically diverging switching rates,
	\item the effect of not amplifying the nominal controller before application, and
	\item instability under slow (but constant) scheduling.% for the nonlinear control system \eqref{eq: SpaceCraft Dyn}. 
\end{itemize}

The attitude dynamics of a rigid spacecraft are modeled by Euler's equations along with a suitable representation of the attitude \cite{ref:astolfi-02}. The attitude dynamics relate the angular velocity of the spacecraft in the body frame to the net external torque acting on the spacecraft, and Euler's equations are derived from the rotational dynamics via Newton's laws. The attitude kinematics, with quaternions representing the attitude, and Euler's equations together are given by
\begin{equation}
\label{eq: SpaceCraft Dyn}
\begin{aligned}
\dot{\quat[]}_0(t) &= -\frac{1}{2} \quat[v]^\top(t) \angM(t) \\ 
\dot{\quat[]}_v(t) &= \frac{1}{2}\big(\quat[0](t)\angM(t) + \quat[v](t) \times \angM(t)\big)\\
\iner\dot{\angM}(t) &= -\angM(t) \times \iner\angM(t) + \conT(t),
\end{aligned}
\end{equation}
where \(\angM(t) \in \R^3 \) is the spacecraft angular velocity expressed in body frame at time \(t\), \( \iner \in 
\R^{3\times3}\) is the (constant) inertia matrix of the rigid body, \(\conT(t) \) is the external torque acting on the spacecraft at time \(t\), and the pair $\bigl(\quat[0](t),\quat[v](t)\bigr)\in\R\times\R^3$ constitutes the unit quaternion at time \(t\), i.e., \(\quat[0]^2(t)+\quat[v]^\top(t)\quat[v](t)=1\). Here the symbol \(\times\) stands for the standard cross product of vectors in \(\R^3\). The external torque \(\conT \) in \eqref{eq: SpaceCraft Dyn} is the feedback control acting on the rigid spacecraft that we design for attitude stabilization. For the attitude stabilization problem the equilibrium point of interest is \(\xEq \Let (\bar{\quat[]},0,0,0) \) with \(\bar{\quat[]} = (1\ 0\ 0\ 0)\). In an obvious way we shift the coordinate system appropriate such that \(0 \in \R^7 \) becomes the closed-loop equilibrium in order to conform with our results. 

In \cite{ref:wen-99} the authors propose an asymptotically stabilizing controller for the spacecraft attitude stabilization problem under the assumption that all control inputs are active and full state information is available at each time; we lift this feedback, given by,
\begin{align}
	\label{eq: Controller}
	\conT(t) = - \const[1]\quat[v](t) - \const[2]\angM(t) \quad\text{for}\ t \geq 0,
\end{align}
for suitable \(\const[1], \const[2] > 0 \), and employ it as our nominal controller. Of course, \(\conT(t) \in \R^3 \) at each \(t\). We employ the feedback \eqref{eq: Controller}, and carry out simulations with the numerical values of different parameters presented in Table \ref{tab: Qtt2} below.
\begin{table}[h]
	\begin{center}
		\begin{tabular}{ @{}ll@{} }
			\toprule
			Quantity & Values \\
			\midrule
			Time of simulation & 200 \si{\second} \\  
			Inertia of spacecraft & \(\textsf{diag}\big((100, 70, 150)\big)\)  \si{\kilogram \meter^2}\\
			\const[1] & 0.5 \\
			\const[2] & 0.1\\
			ODE solver & \(4\)-th order RK \\
			Step length & 0.1 \si{\second}\\
			Initial value of the quaternion & \((1, \ 0,\   0,\  0)^\top\) \\
			Initial value of the angular velocities & \((0.01,\ 0.05,\ 0.03) \) \si{\radian\per \second}\\
			\bottomrule
		\end{tabular}
	\end{center}
	\caption{Parameters for the control system \eqref{eq: SpaceCraft Dyn}; the corresponding controller is presented in \eqref{eq: Controller}.}
	\label{tab: Qtt2}
\end{table}

The switching time is decreased monotonically by a factor of \(0.1 \) starting from \(0.1\si{\second}\) in intervals of \(5 \si{\second} \). Of course, \emph{any} \(\tau(\cdot)\) that is positive, non-increasing and bounded below by the switching signal presented above, would serve the purpose. Asymptotic stability of the equilibrium is guaranteed by Theorem \ref{thm: AS_stable}, and is evident from Fig.\ \ref{fig: quat} and Fig.\ \ref{fig: omga}. 
%The prominent fluctuations in the trajectories of the system discussed in \S\ref{sec: LinModel} are not so pronounced here (cf.\ Fig.\ \ref{fig: omga}) due in part to the (time varying) switching time being smaller than the one used in \S\ref{sec: LinModel}. 

%In the example presented in \S\ref{sec: LinModel}, the removal of the amplification factor of \(m\) (the number of control inputs) cause the trajectories to exponential divergence as discussed in \S\ref{sec: LinModel} and  illustrated in Fig.\ \ref{fig: blowingtraj}. In the present system the effect is not so severe, but leads to a performance deterioration evidenced, e.g., in terms of settling time, as seen in Fig.\ \ref{fig: effectOfm}.  

If the switching frequency is not sufficiently large, then the equilibrium \(\xEq \) can become unstable, and the mechanism by which this happens is fairly subtle. It is well known \cite[Chapter 5]{vidyasagar2002nonlinear} that for an unforced rigid body system for which the intermediate axis is oriented in, say the \(x\) direction, any equilibrium of the form \((\bar{\omega},0,0) \) with \(\bar{\omega} \neq 0\) is \emph{unstable}. For \(\bar{\omega} = 1 \), we tune the controller \eqref{eq: Controller} so that the closed-loop system has an asymptotically stable equilibrium at \(\xEq \Let (\bar{q}, \bar{\omega}, 0, 0) \), which is shifted to \(0 \in \R^7 \) by employing an appropriate shift of the coordinates. Then we periodically schedule the control inputs \eqref{eq: Controller} in a round-robin fashion with a constant switching time of \(18 \si{\second} \). The results are shown in Fig.\ \ref{fig: switchslow}. In this experiment we pick one point each from \(\bddSet_1 \Let \{\bar q\}\times\ball{10^{-1}}{\xEq} \), \(\bddSet_2\Let\{\bar q\}\times\ball{10^{-3}}{\xEq} \) and \(\bddSet_3\Let\{\bar q\}\times\ball{10^{-5}}{\xEq} \) uniformly randomly to avoid any bias. Using these randomly selected points as initial conditions for the closed-loop system, we generate the trajectories in Fig.\ \ref{fig: switchslow}, where for each \(\bddSet_i \) , \(i \in \{1,2,3 \} \), we label the corresponding angular velocity trajectory by \(\omega_0^{(i)}\). Observe that \emph{all} the corresponding trajectories digress away from the ball of radius \(0.1\) around the equilibrium. This experiment, while not a conclusive test of instability of the equilibrium, points to the fact that the aforementioned (fixed) ball around the equilibrium is being transgressed by the trajectories starting inside balls of progressively smaller (by several orders of magnitude) sizes, hinting strongly at instability of the equilibrium under slow switching. Fig.\ \ref{fig: switchslow} also demonstrates that trajectories starting closer to the equilibrium take longer to exhibit divergence away from it. Thus, in order to stabilize the system one needs to reduce the switching time (equivalently, increase the switching frequency) appropriately as dictated by Theorem \ref{thm: MainResult}.

\begin{figure}
	\centering
	\includegraphics[width = 0.6\textwidth,keepaspectratio]{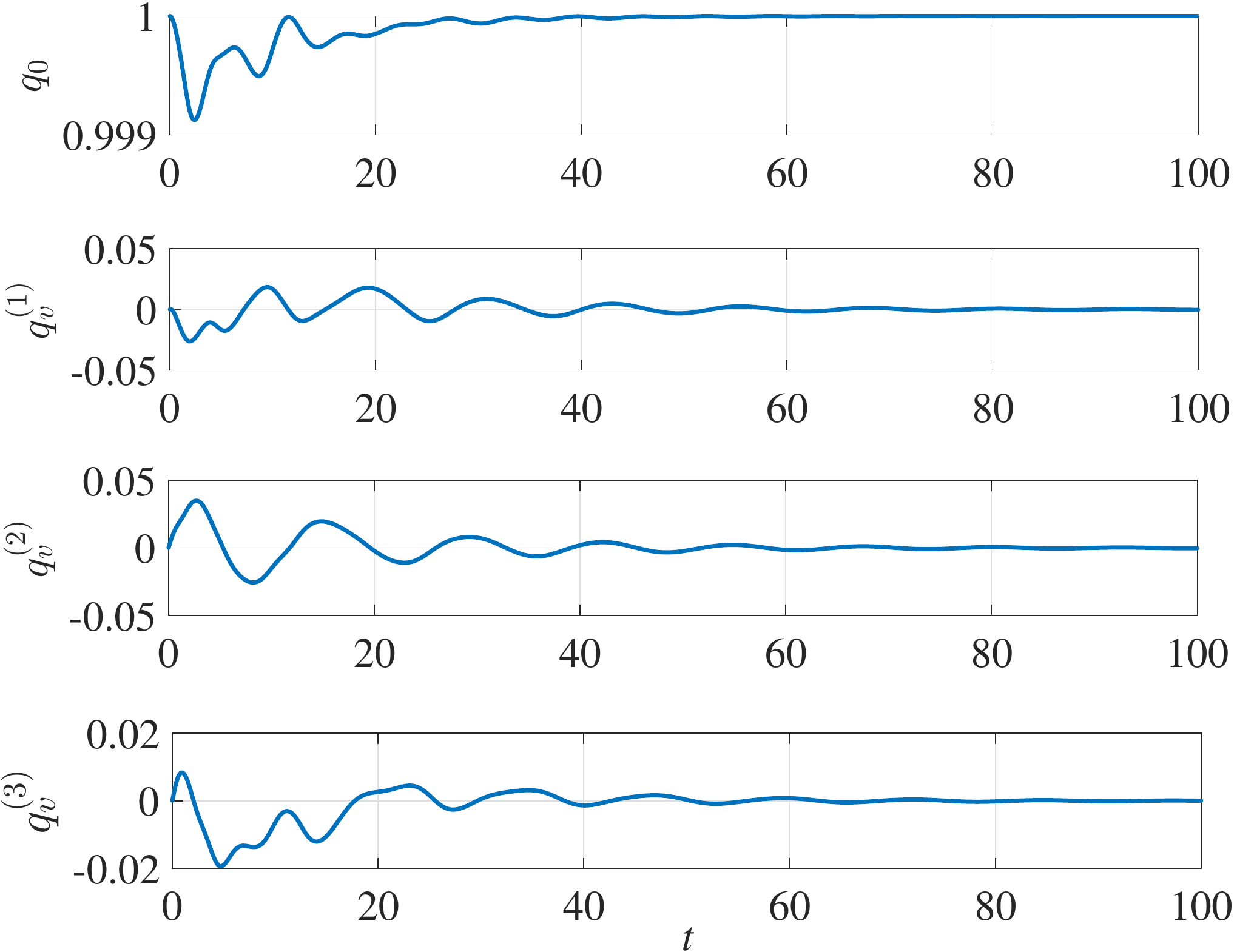}
	\caption{The time evolution of the quaternions; the top subplot corresponds to \(q_0\) and the remaining three correspond to the components of \( q_v\).}
	\label{fig: quat}
\end{figure}
\begin{figure}
	\centering
	\includegraphics[width = 0.6\textwidth,keepaspectratio]{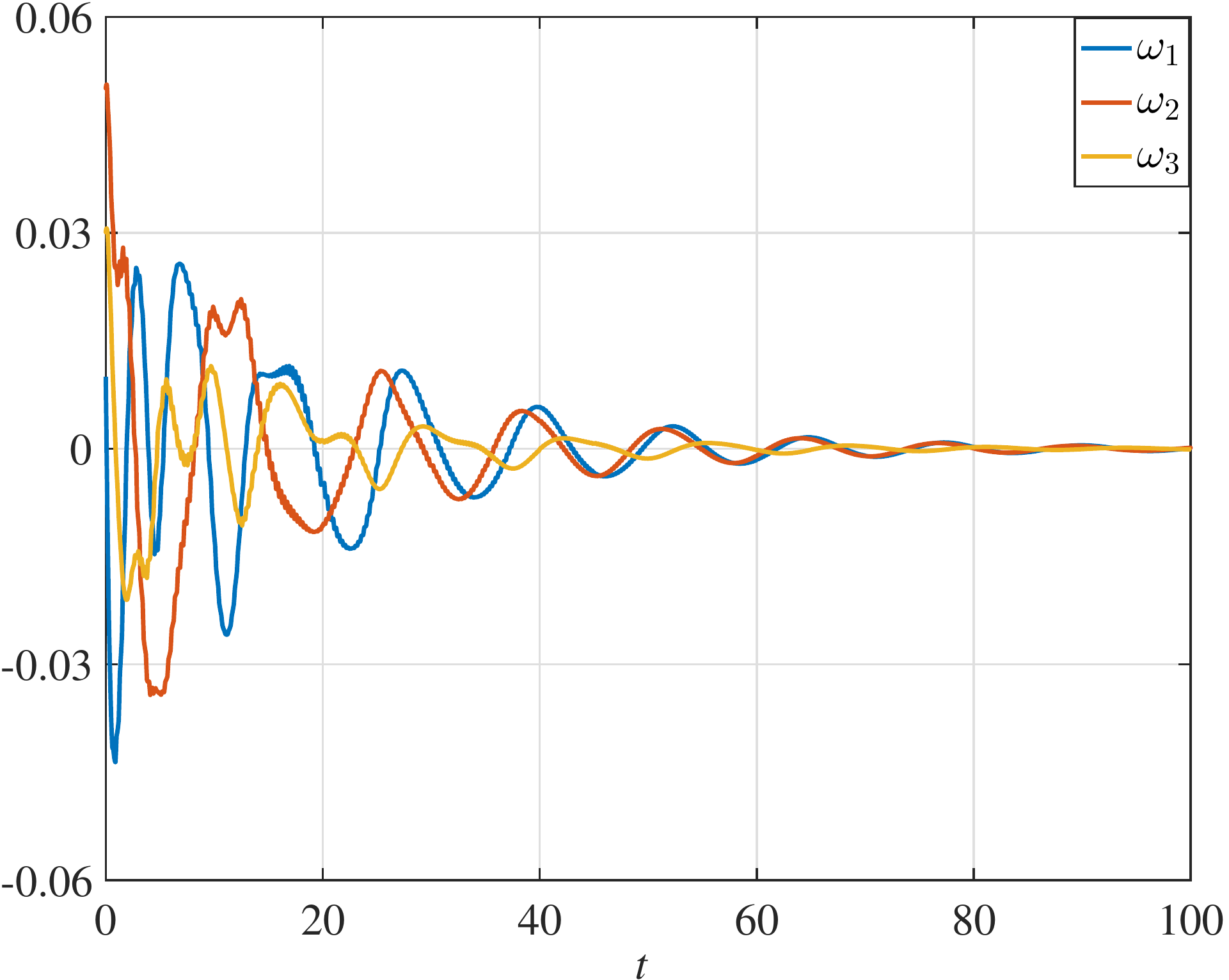}
	\caption{The time evolution of the angular velocities of the spacecraft.}
	\label{fig: omga}
\end{figure}
%\begin{figure}[h!]
%	\centering
%	\includegraphics[width = 0.4\textwidth,keepaspectratio]{1Cont}
%	\caption{The control law implemented based on the periodic scheduling specified in \S\ref{sec:MRS}. The zoomed in sub-plot clearly illustrates the round-robin property of the controller. }
%	\label{fig: ControlSpacecraft}
%\end{figure}   
%\begin{figure}
%	\centering
%	\includegraphics[width = 0.4\textwidth,keepaspectratio]{1effectOfm}
%	\caption{There is a significant qualitative improvement in performance when the control action is amplified. The trajectories plotted in the figure is that of \(\omega_1 \) for both the amplified and the unamplified cases.}
%	\label{fig: effectOfm}
%\end{figure}
\begin{figure}
	\centering
	\includegraphics[width = 0.4\textwidth]{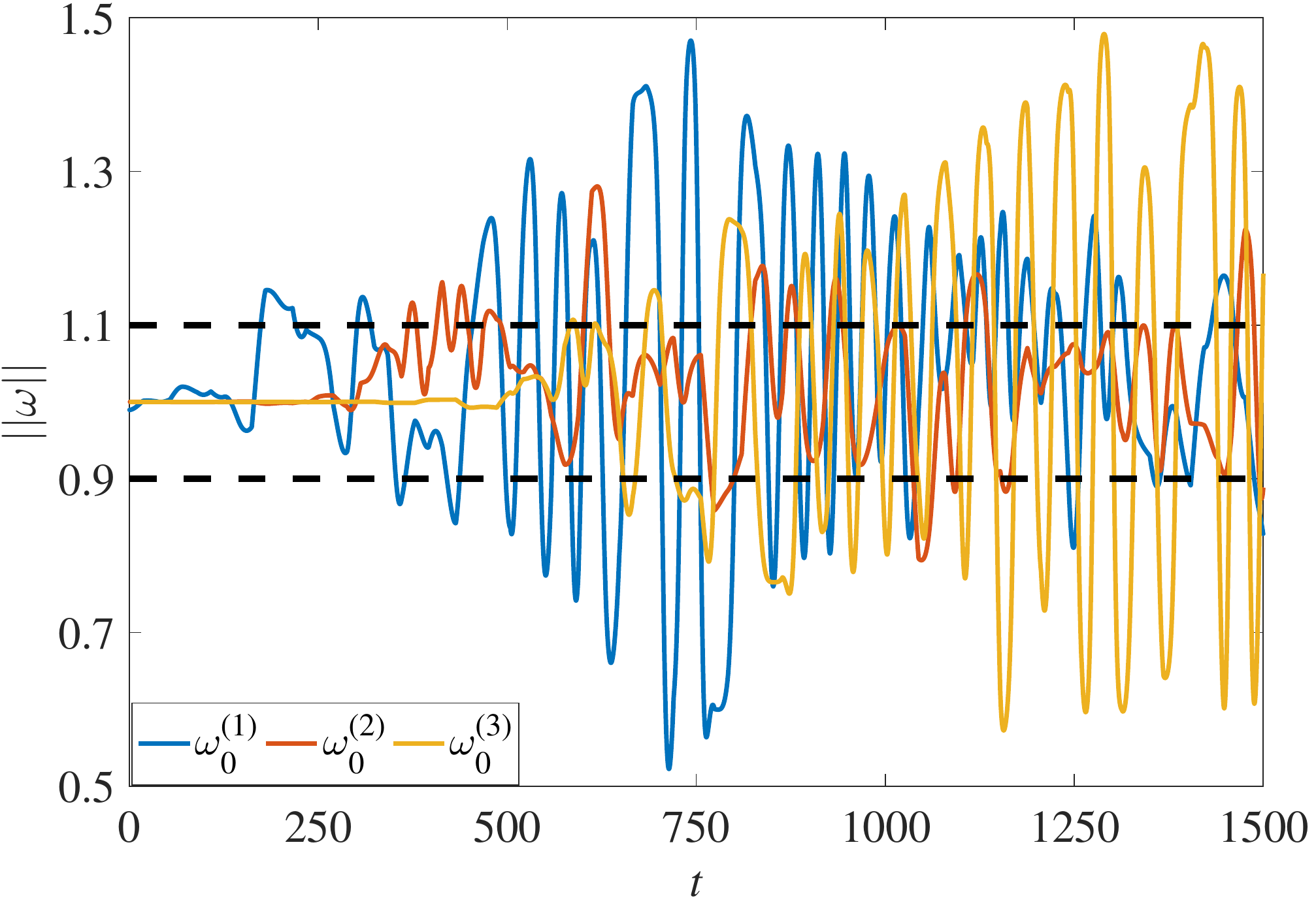}
	\caption{Round-robin scheduling with a constant switching time of \(18\si{\second}\). The trajectories of the spacecraft corresponding to \(\norm{\omega}\) digress away from the \(0.1\)-ball around the equilibrium \( \xEq = (\bar{q}, 1,0,0)\) even if initialized very close to the equilibrium. \(\omega^{(3)}_0\) is the initial condition closest to \(\xEq\) among the three plotted here, and takes the longest to exhibit divergence.}
	\label{fig: switchslow}
\end{figure}

\section{Concluding remarks and future research directions}
\label{sec: Conc}
%%%%%%%%%%%%%%%%%%%%%%%%%%%%%%%%%%%%%%%%%%%%%%%%%%%%%%%%%%%%%%%%%%%%%%%%%%%%%%%%
We proposed two periodic scheduling schemes for a wide class of control-affine systems for which locally asymptotically stabilizing controllers are available, one leading to stability and the other to asymptotic stability. Numerical results were provided to illustrate the results. 

A natural line of investigation is to quantify the performance of controllers and their convergence with respect to the frequency of the periodic scheduling. Moreover, it is of significant interest to study the behavior of such systems under randomization as opposed to the round-robin scheduling treated herein. These topics will be taken up in subsequent works.
\appendix
%%%%%%%%%%%%%%%%%%%%%%%%%%%%%%%%%%%%%%%%%%%%%%%%%%%%%%%%%%%%%%%%%%%%%%%%%%%%%%%%
\section{Essential Ingredients}
\label{app: PTheorem}
%%%%%%%%%%%%%%%%%%%%%%%%%%%%%%%%%%%%%%%%%%%%%%%%%%%%%%%%%%%%%%%%%%%%%%%%%%%%%%%%
This appendix collects some standard and not-so-standard results that are central to the proofs of Theorem \ref{thm: MainResult},  \ref{thm: AS_stable},  \ref{thm: NonuniformResult}, and \ref{thm: RegionOfAttractionUniformResult} presented in \S\ref{app: Proof}. 

%%==============================================================================
\subsection{Chattering lemma}
%%==============================================================================
We start with a version of the so-called ``chattering lemma'' \cite[Theorem 3.6.1]{ref:BerMed-13}, which is popular in optimal control theory, applied to time-invariant right-hand sides.% (the sense of ``time'' will become clear in the following theorem)
\begin{theorem}\label{thm: ChatterLemma}
	Let \(n,q \in \N\), let \(\compSet \subset \R^n \) be a compact and nonempty set and let \(\interval \subset \R \) be a compact interval. Define the maps \(\compSet \ni x \mapsto \ithVecField(x) \in \R^n, \) \(i \in \setInt[\numSys]\), satisfying the following properties:
	\begin{enumerate}[leftmargin=*, label = {(C-\roman*)}, align = left, widest = ii]
		\item each \(\ithVecField\) is continuous in \(\compSet\);
		\item there exist constants \(\lipConst_1, \lipConst_2 > 0\)  such that for all \(x, y\in \compSet\) and \(i \in \setInt[\numSys] \):
		\begin{align*}
			\norm{\ithVecField(x)} &\leq \lipConst_1, \\
			\norm{\ithVecField(x)-\ithVecField(y)} &\leq \lipConst_2 \norm{x-y}.
		\end{align*}
	\end{enumerate}
	Let \(\interval \ni t \mapsto \partUnity(t)\in\lcro{0, +\infty}\), \(i \in \setInt[\numSys] \), be measurable maps satisfying \(
	\sumFunc{i}{\numSys}\partUnity[i](t) = 1
	\) almost everywhere.
	Then for every \(\epsBar\ > 0\) there exists a subdivision of \(\interval\)  into a finite collection \(\{\lenInt[j]\}_{j=1}^{\numDiv}\) of non-overlapping intervals and an assignment of one of the functions \(\ithVecField[1], \ithVecField[2], \dots, \ithVecField[\numSys]\) to each \(\lenInt[j]\)
	%	 such that  if \(\ithVecField[\lenInt]\) denotes the function assigned to the interval \(\lenInt\) and if the map \(\interval \times \compSet \ni (t,x) \mapsto \ithVecField[](t,x) \in \R^n\) is defined such that %agree with \(\ithVecField[\lenInt]\) on the interior \(\inLenInt[\lenInt]\) of each \(\lenInt\), that is, 
	%	\begin{align*}
	%	\ithVecField[](t,x) = \ithVecField[\lenInt](x) \quad \text{for} \ t \in \inLenInt[\lenInt], x \in \compSet , \text{and}\ j \in \setInt[\numDiv],
	%	\end{align*}
	then for every \(t', t''\) in \(\interval\) and all \(x \in \compSet\) 
	\begin{align}
		\norm[\bigg]{\int_{t'}^{t''}\Big( \sumFunc{i}{\numSys} \partUnity[i] (t) \ithVecField[i](x) - \ithVecField[\sigma(t)](x)\Big) dt} < \epsBar,
	\end{align} 
	where the assignment is characterized by the map \(\interval \ni t \mapsto \sigma(t) \in \setInt[q] \).
\end{theorem}

We refer to the variable \(t\) as the \emph{time} argument and the variable \(x\) as the \emph{space} argument in the sequel.%In lieu of this we call the functions \(\ithVecField[i]\) to be time invariant \((\)for every \(i \in\setInt[\numSys])\) and the function \(\ithVecField[] \) to be time dependent. 

\begin{remark}\label{rem: Chatterassignment}
	{\rm The candidate subdivision of the compact interval \(\interval\) and an assignment of the functions \(\ithVecField[i]\), for \(i \in \setInt[\numSys] \) into those subdivisions, whose existence in stated in Theorem \ref{thm: ChatterLemma}, as used in \cite{ref:BerMed-13} is presented here. Define \(\eps \Let \frac{\epsBar}{2(2+|\interval|)}\). For an integer \(\tilde{\numDiv}\), we define a partition of \(\interval \) into \(\big\{\lenIntC\big\}_{j \in \setInt[\tilde{\numDiv}]}\), such that \(\abs{\lenIntC} < \frac{\eps}{\max\{\lipConst_1, \lipConst_2\}}\). Further, for each \(j \) the interval \(\lenIntC\) is subdivided into finitely many non-overlapping sub-intervals \(\subDiv{j}{1}, \subDiv{j}{2}, \dots, \subDiv{j}{\numSys}\) such that 
		\begin{align*}
			\abs{\subDiv{j}{i}} \defas \int_{\lenInt[j]}^{} \partUnity(t) dt \quad \text{for}\ i \in \setInt[\numSys].
		\end{align*}
		The assignment follows the rule 
		\begin{align*}
			\sigma(t) = i \quad \text{for}\ t \in \inLenInt[\subDiv{j}{i}], \ \ i \in \setInt[\numSys], \   j \in \setInt[\numDiv]. 
		\end{align*}
		We shall be employing this particular subdivision in the sequel.}
\end{remark}	

%%==============================================================================
\subsection{Continuous dependence of solutions of a differential equation on parameters}
%%==============================================================================
We state a result on the dependence of the solutions of a parameter-dependent differential equation on the parameter. For a positive scalar \(\initT \) and some positive integers \(d \) and \(p\), consider the two differential equations
\begin{align}\label{eq: filSys}
	\dot{\filStateX}(t) = \filDE{\filStateX(t)}{\filParam}{t},  \quad  \filStateX(\initT) = \stinit, \quad t \geq \initT,
\end{align}
\begin{align}\label{eq: filSys2}
	\dot{\filStateY}(t) = \filDEY{t}{\filStateY(t)},  \quad  \filStateY(\initT) = \stinit,\quad t \geq \initT,
\end{align}
where \(\filStateX(t), \filStateY(t) \in \R^d\) for each \(t\), the (fixed) parameter \(\filParam\) takes values in \(\filParamSet\), and the maps \(\lcro{t_0, +\infty} \times \R^d \times \R^p \ni (t,x,\filParam) \mapsto \filDE{\filState}{\filParam}{t}  \in \R^d \) and \(\lcro{t_0, +\infty} \times \R^d \ni (t,x) \mapsto \filDEY{t}{x}  \in \R^d \)  satisfy the standard Carath\'eodory conditions \cite[Chapter 1]{filippov2013differential}.
%This result provides a set of sufficient conditions that ensures continuous dependence of solution trajectories of differential equation \eqref{eq: YMeasForm} on a parameter. One set of sufficient conditions that the solution trajectories must satisfy for such continuous dependence on some parameter is presented in Theorem \ref{thm: filLemma}. However, it becomes difficult to verify these requirements when the solution trajectories of a differential equations cannot be computed or is difficult to compute (which is often the case with nonlinear systems). This difficulty is circumvented in the Theorem \ref{thm: FilTheorem}, where one requires to only check certain conditions on the right-hand side of the differential equation to guarantee continuous dependence of solution trajectories on certain parameter.
%For positive integers \(d\) and \(\filPPS\), consider the differential equation driven by a Young measure \(\YMeasP[\filParam] \), which depends on some parameter \(\filParam \in \filParamSet \), given by
%\begin{align}\label{eq: YMeasForm1}
%\dot{\filState}(t) =  \avgof{{\YMeasP[\filParam]}}(t,\filState(t)), \quad\filState(\timeVal[0]) = \stinit,
%\end{align}
%for \(t \geq \timeVal[0] \). 
We assume that both \eqref{eq: filSys} and \eqref{eq: filSys2} have unique solutions and we denote these solutions by \(\lcro{t_0, +\infty} \times \R^p \ni (t, \filParam) \mapsto \filStateX(t,\filParam) \in \R^d \) for \eqref{eq: filSys} and by \(\lcro{t_0, +\infty} \ni t \mapsto \filStateY(t) \in \R^d \) for \eqref{eq: filSys2} , where we suppress the dependence of solution trajectories on both the initial time instant \( \initT\) and the initial condition \(\stinit \) for the sake of brevity. The following theorem provides a set of sufficient conditions for the convergence of the trajectories of \eqref{eq: filSys} to those of \eqref{eq: filSys2} on a compact interval.  
%, the Young measure \(\YMeasP[\filParam_i] \xrightarrow[i \sra +\infty]{\text{weak}-\star} \YMeasF\).
% We represent the solution trajectories of the differential equation governed by the average of the Young measure \(\YMeasF\) as \(\R \ni t \mapsto \filState_{\eP}(t) \in \R^d\).
\begin{theorem}\cite[Chapter 1, Theorem 7]{filippov2013differential}\label{thm: FilTheorem}
	Consider the differential equations \eqref{eq: filSys} and \eqref{eq: filSys2}. Define a compact interval \(\interval \Let [t_0, t_1] \) for some scalars \(t_0 < t_1 \) and a compact set \(\compFil \subset \R^d\) containing the point \(\stinit \). Let
	\begin{enumerate}[leftmargin=*, align= left, label = (A-\roman*),widest = iii]
		\item \label{th: avgFirst}the function \( \filDE{\filState}{\filParam}{t}\) be measurable in \(t\) for constant \(x, \filParam \); and the function \( \filDEY{t}{\filState} \) be measurable in \(t\) for constant \(x\);
		
		\item \label{th: avgSecond}\(\norm{\filDE{\filState}{\filParam}{t}} \leq \filUB\), where \(\R \times \R^\filPPS \ni (t, \filParam) \mapsto \filUB \in \R \) is summable in \(t\); similarly \(\norm{\filDEY{t}{\filState}} \leq \filUBF\), where \(\R \ni t \mapsto \filUB \in \R \) is summable in \(t\);
		
		\item \label{th: avgThird}there exists a summable function \(\interval \in t \mapsto \filTimeCap \in \R \) and a  monotone function \(\R^d \ni x \mapsto \filSpaceCap \in \R \) satisfying $\ \filSpaceCap\sra 0 \ \text{as} \ \filState \sra 0$, such that for each \(r > 0\), for all \(x, y \in \compFil \) satisfying \( \norm{x-y} \leq r \) and for almost all \(t \in \interval\) we have
		\begin{equation}
		\begin{aligned}\label{eq: LipNormfil}
		\norm{\filDE{x}{\filParam}{t} - {\filDE{y}{\filParam}{t}}} &\leq \filTimeCap \filSpaceCap[r] \ \text{and} \\ 
		\norm{\filDEY{t}{x} - {\filDEY{t}{y}}} &\leq \filTimeCap \filSpaceCap[r];
		\end{aligned} 
		\end{equation}
		\item \label{th: avgFourth}for each \(\filState \in \compFil\) and any sequence \(\filParam_i \xrightarrow[i \sra +\infty]{} \filParam_{\eP}\), we have
		\begin{align}\label{eq: 19}
			\int_{t_0}^{t} \filDE{x}{\filParam}{s}\ \dd s \xrightarrow[i \sra +\infty]{} \int_{t_0}^{t} \filDEY{s}{x}\ \dd s
		\end{align}
		uniformly in \(t\) on the interval \(\interval\).
	\end{enumerate}
	%	Then for any \(\filParam \) sufficiently near \(\filParam_{\eP}\), the solution of the differential equation \eqref{eq: YMeasForm} on the interval \(\interval\) exists and converges uniformly to \(\filState_\eP(\cdot)\) on \(\interval \) as  \(\filParam= \filParam_i \overset{i \sra +\infty}{\ra} \filParam_{\eP} \).
	Then for every \(\eps > 0 \) there exists a scalar \( \etta >  0\) such that for all \(\filParam\) satisfying \(\norm{\filParam - \filParam_{\eP}} < \etta \), each solution \(\filStateX(\cdot,\filParam)\) of \eqref{eq: filSys} \((\)corresponding to the initial condition \(\stinit \in \compFil)\) exists on the interval \(\lcrc{t_0, t_1} \) and satisfies 
	\[
	\norm{\filStateX(t,\filParam)- \filStateY(t)} < \eps \quad \text{for}\ t \in \lcrc{t_0, t_1}\ \text{and}\ \norm{\filParam - \filParam_{\eP}} < \etta.
	\]
\end{theorem}

\subsection{A version of Alekseev's bound}
%%==============================================================================
At this juncture we state a result from \cite{ConcentrationBound} that provides a bound on the difference between solution trajectories of a class of autonomous differential equations when initialized	 at two different initial conditions. For a positive integer \(d\) and a scalar \(t_0\), consider the following differential equation
\begin{align}\label{eq: sysAsymp}
	\dot{x}(t) = \sysAsymp(x(t)), \quad x(\initT) = \stinit, \quad t \geq \initT, 
\end{align}
where the map \(\R^d \ni x \mapsto \sysAsymp(x) \in \R^d\) is twice continuously differentiable.
%\begin{enumerate}[leftmargin = *, label = (\roman*), widest = ii, align = left]
%\item \label{ass: Thoppe1} the map \(\R^d \ni x \mapsto \sysAsymp(x) \in \R^d\) is continuously differentiable.
%\end{enumerate}
We denote the solution of \eqref{eq: sysAsymp} by \(t\mapsto x(t,t_0, \stinit) \) for \(t\geq t_0\). Let \(0 \in \R^d\) be an \emph{isolated, hyperbolic} and  \emph{locally asymptotically stable} equilibrium of \eqref{eq: sysAsymp} and let \(\compTho\) be a bounded set containing \(0\) in its interior and that is contained in the domain of attraction of the equilibrium point \(0\). 

Recall that a continuously differentiable function \(\lyap: \domainConv \subset \R^d \ra \R\) is said to be a \emph{Lyapunov function} for the equilibrium point \(0\) if \(\domainConv\) is open, \(\lyap(0) = 0\), and for all \( x \neq 0 \), the estimates \(\lyap(x) > 0\) and \( \inprod{\nabla\lyap(x)}{\sysAsymp(x)} < 0 \) hold simultaneously. The converse Lyapunov Theorem \cite[Chapter 5]{vidyasagar2002nonlinear} guarantees the existence of a Lyapunov function near a locally asymptotically stable equilibrium point \(0\). In fact, we can choose \(\lyap\) so that $\lyap(x) \to +\infty$ as \(x\) tends to the boundary of the \(\domainConv\). Let \(\lyap\) be a Lyapunov function for \eqref{eq: sysAsymp} such that \(\domainConv\) denotes the domain of the Lyapunov function \(\lyap\). For \(p > 0 \) we define \(\lyapBall[p] \Let V^{-1}(p) \), and for \(q>0\) we let \(	\Dilate{q}{\lyapBall[p]} \Let \bigcup_{y \in \lyapBall[p]} \ball{q}{y}.\)
From elementary properties of Lyapunov functions we see that there exist scalars \(r, r_0, \eps_0\) such that \(r > r_0\) and for each \(\eps \in \loro{0, \eps_0} \) we have
\begin{align}
	\label{eq: rDef2}
	\ball{\eps}{0} \subset \compTho \subset \lyapBall[r_0] \subsetneq \Dilate{\eps_0}{\lyapBall[r_0]} \subset \lyapBall[r] \subsetneq \domainConv.
\end{align}
A pictorial representation of the above statement is presented in Fig.\ \ref{fig: picSet}.
\begin{center}
	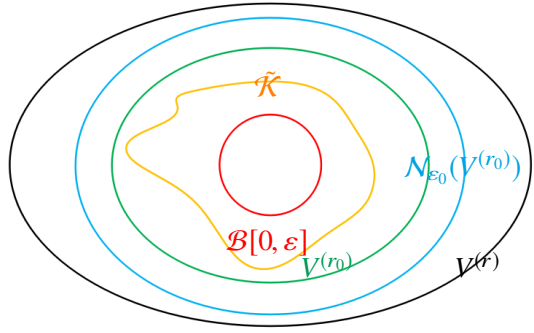
\begin{figure}
		\begin{tikzpicture}
		\node[anchor=south west,inner sep=0] (image) at (0,0) {\includegraphics[width=0.6\textwidth]{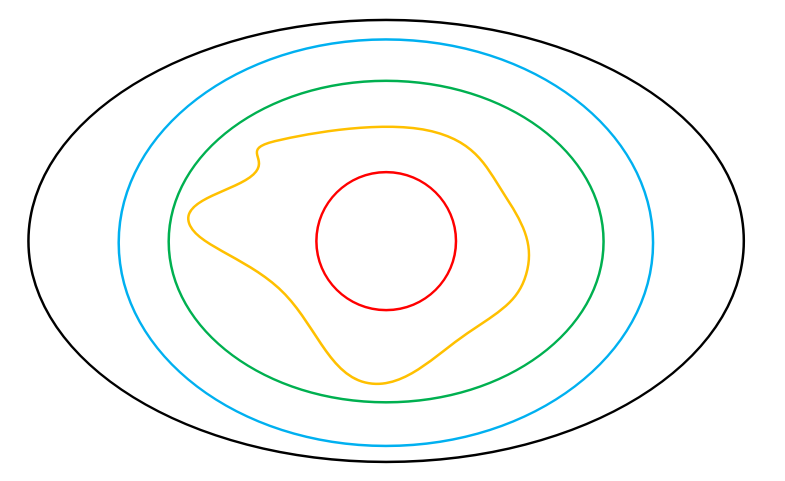}};
		%	\begin{scope}[x={(image.south east)},y={(image.north west)}]
		%	\draw[help lines,xstep=.1,ystep=.1] (0,0) grid (1,1);
		%	\foreach \x in {0,1,...,9} { \node [anchor=north] at (\x/10,0) {0.\x}; }
		%	\foreach \y in {0,1,...,9} { \node [anchor=east] at (0,\y/10) {0.\y}; }
		%	\end{scope}
		\node at (3.7,1.3) {{\color{red}\(\ball{\eps}{0}\)}};
		\node at (3.7,3.4) {{\color{orange}\(\compTho\)}};
		\node at (4.5,1) {{\color{ForestGreen} \(V^{(r_0)} \)}};
		\node at (6.3,2.3) {{\color{Cerulean}\(\Dilate{\eps_0}{\lyapBall[r_0]} \)}};
		\node at (6.5,1) {\(\lyapBall[r] \)};
		\end{tikzpicture}
		\caption{Pictorial depiction of the sets in \eqref{eq: rDef2}.}
		\label{fig: picSet}
	\end{figure}
\end{center}

Let \(\jacobian[0]\) be the Jacobian matrix of \(\sysAsymp\) evaluated at \(0\). Since the equilibrium point \(0\) is hyperbolic and locally asymptotically stable, the matrix \(\jacobian[0]\) is guaranteed to be Hurwitz.\footnote{The assumption of the equilibrium point being \(0\) is missing from \cite{ConcentrationBound}.} Let \(\eigVal[1], \eigVal[2], \dots, \eigVal[d]\), with possible repetitions, denote the eigenvalues of the matrix \(\jacobian[0]\). We define \(\minEig \Let \min_{i\in \setInt[d]} \big\{ -\Re(\eigVal)\big\}. \)
Fix two scalars \(\lamPrime, \kap > 0\) such that \(\lamPrime < \minEig \) and \(0 < \kap < 1 \). Then \cite[Corollary 3.6]{teschl2012ordinary} asserts that there exists a constant \(\constExp\ >0 \) such that 
\(
\bigl\|{\expo[{\jacobian[0] t}]}\bigr\| \leq \constExp \expo[-\lamPrime t] \ \text{for all}\ t \geq 0.
\)
\begin{theorem}{\cite[Lemma 2]{ConcentrationBound}}
	\label{thm: thoppeResult}
	Consider the system \eqref{eq: sysAsymp} along with its associated data and adopt the notations developed above. Define
	\begin{equation}\label{eq: LamDef}
	\lam \Let \frac{1-\kap}{\constExp^2} \lamPrime.
	\end{equation}
	Let \(r\) be chosen as stated in \eqref{eq: rDef2}, let \(\initCond[1], \initCond[2]\) be arbitrary points in \(\lyapBall[r]\), and let \(s \ge 0\). Then
	\begin{align*}
		\norm{x(t,s,\initCond[1]) - x(t,s,\initCond[2])} \leq \constLemma \norm{\initCond[1] - \initCond[2]} \expo[-\lam(t-s)] \quad \text{for all}\ t\geq s,
	\end{align*}
	where \(\constLemma\) is a constant that depends on the set \(\compTho\) and the vector field \(\sysAsymp\).
\end{theorem}
%\begin{remark}
%	If one chooses larger \(\kap\) then the constant \(\constLemma\) decreases and we get a stricter bound. Reader is referred to \cite[Section 5]{ConcentrationBound} for further reading. 
%\end{remark}

%%%%%%%%%%%%%%%%%%%%%%%%%%%%%%%%%%%%%%%%%%%%%%%%%%%%%%%%%%%%%%%%%%%%%%%%%%%%%%%%
\section{Proofs of the Main Results}
\label{app: Proof}
%%%%%%%%%%%%%%%%%%%%%%%%%%%%%%%%%%%%%%%%%%%%%%%%%%%%%%%%%%%%%%%%%%%%%%%%%%%%%%%%

\begin{tcolorbox}
	For a consistent notation for next two subsections, we shall denote the trajectories of  \eqref{eq: normalSys} by \(\stateTA(\cdot) \), of \eqref{eq: sparseSys} by \(\stateSS(\cdot) \), and of  \eqref{eq: sparseSysAS} by \(\stateAS(\cdot)\).
\end{tcolorbox}

\subsection{Proof of Theorem \ref{thm: MainResult}}\label{app_lab: Main Result}
The proof unfolds in the following two key steps:
\begin{enumerate}[leftmargin=*, align=left, label = (M-\roman*),widest = ii]
	\item\label{step:1} We bound the difference between trajectories of \eqref{eq: normalSys} and \eqref{eq: sparseSys} on a compact time interval. To accomplish this, we appeal to Theorem \ref{thm: FilTheorem}; 
	\item\label{step:2} We iteratively extend the preceding bound over the entire time horizon \( \lcro{\initT, +\infty}\) by repeatedly invoking Theorem \ref{thm: FilTheorem} in conjunction with Theorem \ref{thm: thoppeResult}. Local asymptotic stability of the equilibrium \(0 \in \R^d \) of the closed-loop system \eqref{eq: normalSys} plays a crucial role in the arguments of the proof.  
\end{enumerate}
We proceed to the comprehensive proof of Theorem \ref{thm: MainResult}.

\noindent\(\bullet\) \textbf{Step-I:} As the equilibrium \(0 \in \R^d \) of \eqref{eq: normalSys} is asymptotically stable, we can find a \(\initAsymp > 0 \) such that if \(\normInit \leq \initAsymp \) then \(\norm{\stateTA(t)} \xrightarrow[t \rightarrow +\infty]{} 0\). Fix \(\ups > 0 \), and define \(\betBound \Let \min(\ups, \initAsymp) \). Define a compact set  \(\bddSet \Let   \ball{\betBound}{0} \). Pick a scalar \(\bound \in \loro{0,1} \) and define the another compact set \(\bddSet_\bound \Let \ball{\bound\betBound}{0}\). Also, the equilibrium point \(0 \in \R^d \) of \eqref{eq: normalSys}  is stable in the sense of Lyapunov therefore we can find a positive scalar \(\del \) such that if \(\normInit \leq \del \) then \(\stateTA(t) \in \bddSet_{\alpha} \) for all \(t \geq \initT \): in other words, the trajectory \(\stateTA(\cdot) \subset \bddSet_{\alpha} \subsetneq \bddSet \).
%As per the assumptions stated in Theorem \ref{thm: MainResult}, the equilibrium \xEq\ of the closed-loop system \eqref{eq: normalSys} is locally asymptotically stable. By fixing a scalar \(0 < \scalarProp < 1 \) as in Proposition \ref{prop: ClaimMy} for the closed-loop dynamics \eqref{eq: normalSys}, we can find a scalar \(\timeDown[] > 0\) such that after this time the asymptotically stable trajectory, \(x(\cdot) \), is contained in a smaller ball for all time in the future. That is, for the scalars \(q, s > 0\) if \(x(s)\) is contained in  \( \ball{q}{\xEq}\) then for all \(t \geq s+\timeDown[]\), we have \(x(t) \in \ball{\scalarProp q}{\xEq} \). Without loss of generality, we shall take \(\xEq = 0 \) and the scalar \timeDown[]\ to be a multiple of \(\contDim\switchTime \) in the forthcoming proof. That is, \(\timeDown[] = \intFac\cdot \contDim \switchTime \), where \(\intFac\) is  a positive integer.	
%	We denote a Lyapunov function for the closed-loop system \eqref{eq: normalSys} by \(\R^\stateDim \ni x \mapsto \lyap(x) \in \R \).

For \(\timeDown[]>0 \), define \(\FAC \Let \constLemma \expo[{-\lam\timeDown[]}]\), where \(K_1, \lam \) are positive scalars which depend on the set \(\bddSet\) and the closed-loop vector field governing \eqref{eq: normalSys}. Choose the scalar \(T\) such that \(\FAC < 1 \). We partition the interval \([\initT, +\infty[ \) into a countable family of disjoint half-open intervals of length \(\timeDown[]\) as \(\lcro{\initT, +\infty} = \bigcup_{n\in\nat} \intV[n]\) where \( \intV[n] \Let  \lcro{\initT + n \timeDown[], \initT + (n+1)\timeDown[]} \). We now provide a  uniform bound on the difference between the trajectories of \eqref{eq: normalSys} and \eqref{eq: sparseSys} over a compact interval by appealing to Theorem \ref{thm: FilTheorem}. 

The requirements of Theorem \ref{thm: FilTheorem} except \ref{th: avgFourth} are trivially met by the right hand sides of \eqref{eq: normalSys} and \eqref{eq: sparseSys} since these are the standard Carath\'eodory conditions. 
We claim that \ref{th: avgFourth} is also satisfied in view of Theorem \ref{thm: ChatterLemma} and Remark \ref{rem: Chatterassignment}. To verify this we define the map
\begin{align*}
	\R^\stateDim \ni x \mapsto \ithVFMainRes(x) \Let \drift[x] +  \contDim \channel{i}{x} \control[i](x) \in \R^\stateDim,
\end{align*}  
a partition of unity as \(\lcro{\initT, +\infty} \ni t \mapsto \partUnity(t) = \frac{1}{\contDim} \in \lcro{0, +\infty}\) for each \(i \in \setInt[\contDim] \), and the time-varying map 
\begin{align*}
	\lcro{\initT, +\infty} \times \R^d \times \loro{0,+\infty} \ni (t,x, \switchTimeT) \mapsto \VFMainRes(t,x, \switchTimeT) \Let \ithVFMainRes[\switchFunc(t-\initT,\switchTimeT)](x) \in \R^d, 
\end{align*}
where \(\switchFunc(\cdot) \) is as defined in \eqref{eq: SwitchFunc}.   
These definitions permit us to write \eqref{eq: normalSys} as
\begin{align*}
	\dot{\stateTA}(t) = \sumFunc{i}{\contDim} \partUnity[i](t) \ithVFMainRes[i](\stateTA(t)), \quad \stateTA(\initT) = \stinit, \quad t \geq \initT,
\end{align*}
and \eqref{eq: sparseSys} as 
\begin{align*}
	\dot{\stateSS}(t) = \VFMainRes(t,\stateSS(t), \switchTime), \quad \stateSS(\initT) = \stinit, \quad t \geq \initT. 
\end{align*}
Moreover, from Theorem \ref{thm: ChatterLemma} for \(t \in \intV[k] \) (\(k \in \nat \)) and all \(\zSt \in \bddSet\) we have 
\begin{align}\label{eq: chatterK}
	\int_{\initT+k\timeDown[]}^{t} \bigg( \sumFunc{i}{\contDim} \partUnity[i](s) \ithVFMainRes[i](\zSt) - \VFMainRes(s,\zSt,\switchTime)\bigg)\dd s \xrightarrow[\switchTime \downarrow 0]{} 0 \quad \text{uniformly in} \ t,
\end{align}
which confirms the claim.

Against the preceding backdrop, it is immediate from Theorem \ref{thm: FilTheorem}  that the trajectory \(\stateSS(\cdot) \) converges uniformly to \(\stateTA(\cdot) \) on the interval \(\intV[0] \).  That is, for each \(\gamBar > 0 \) there exists a switching time \(\switchTimeS\) sufficiently small such that 
\begin{align}\label{eq: unifBoundMR}
	\sup_{t \in \intV[0]} \norm{\stateTA(t) - \stateSS(t)} \leq \gamBar. 
\end{align}

At this stage pick \(\gamBar < \frac{\betBound(1 - \alpha)}{\consLem} \) with \(\consLem \Let 1 + \frac{\constLemma}{1-\FAC} \), and corresponding to this value of \(\gamBar \) we select \(\switchTimeS\) such that \eqref{eq: unifBoundMR} is satisfied, and keep this \(\switchTimeS\) fixed through the end of the current proof.

For the next step of the proof we shall retain the notation of \(y(\cdot) \) for the solution trajectory of \eqref{eq: sparseSys} but we emphasize that this trajectory corresponds to the switching time \(\switchTimeS \) fixed above.

\noindent\(\bullet \) \textbf{Step-II:} In order to extend the uniform bound in \eqref{eq: unifBoundMR} to \(\lcro{\initT, +\infty} = \bigcup_{k \in \N} \intV[k] \), we appeal to the results of Theorem \ref{thm: FilTheorem} over successive intervals \(\intV[k] \). Note that it was possible to use the results of Theorem \ref{thm: FilTheorem} to uniformly bound the separation between two trajectories \(\stateTA(\cdot) \) and \(\stateSS(\cdot) \) over \(\intV[0] \) because the two trajectories start from the same state at time \(\initT\) (along with other requirements posited in Theorem \ref{thm: FilTheorem}). However, once the system has evolved over time interval \(\intV[0] \),  the two trajectories \(\stateTA(\cdot) \) and \(\stateSS(\cdot) \) may not necessarily intersect at any future time, and thus there is no direct way to use Theorem \ref{thm: FilTheorem}. To circumvent this problem, we construct an ensemble of trajectories \(\big(\nthStateTA(\cdot)  \big)_{n\in\N} \), each of which satisfies the dynamics given in \eqref{eq: normalSys} with the initial conditions
\begin{align}\label{eq: InitCond}
	\nthStateTA[n](\initT+n\timeDown[]) = \stateSS(\initT+n\timeDown[]) \quad\text{for each}\ n \in \N.
\end{align}  
That is, for \(t \in \intV[k] \), the differential equation governing \(\nthStateTA[k](\cdot) \)  is 
\begin{align*}
	\dot{\stateTA}_{(k)}(t) = \sumFunc{i}{\contDim} \partUnity[i](t) \ithVFMainRes[i](\nthStateTA[k](t)), \quad \nthStateTA[k](\initT+k\timeDown[]) = y(\initT+k\timeDown[]), \quad t \in \intV[k].
\end{align*}
%and \eqref{eq: sparseSys} as 
%\begin{align*}
%\dot{\stateSS}(t) = \VFMainRes(t,\stateSS(t), \switchTime), \quad 
%\end{align*}
Employing the fact that the vector fields \(\ithVFMainRes \) are time-invariant along with the vector field \(\VFMainRes \) being invariant under time shift by \(m\switchTimeS \) and \eqref{eq: InitCond}, we can conclude from Theorem \ref{thm: FilTheorem} that
\begin{align}\label{eq: boundIntk}
	\sup_{t \in \intV[k]} \norm{\nthStateTA[k](t) - \stateSS(t)} \leq \gamBar \quad \text{for each}\ k \in \N, 
\end{align}
where the switching time is chosen to be \(\switchTimeS \), which is identical to the one selected in \eqref{eq: unifBoundMR}. 
%	Using Theorem \ref{thm: thoppeResult}, we obtain a bound on the difference between the trajectories \(\stateTA(\cdot) \) and \(\nthStateTA(\cdot)\) in the interval \(\intV[n+1] \). 

We claim that the difference between \(\stateSS(\cdot) \) and \(\stateTA(\cdot) \) at \(t = \initT+n\timeDown[] \) for every \(n \in \N\) satisfies the estimate
\begin{align}\label{eq: nthSum}
	\norm{\stateSS(\initT+n\timeDown[]) - \stateTA(\initT+n\timeDown[])} \leq \gamBar\bigg(\sum_{i=0}^{n-1} \constLemma^i\expo[{-i \lam \timeDown[]}] \bigg).
\end{align}
Recall that \(K_1\) is positive scalar which depend on the set \(\bddSet\) and the closed-loop vector field governing \eqref{eq: normalSys}. 

The preceding claim is established inductively. For the induction base, consider \(n =2 \). At \(t = \initT+ 2 \timeDown[] \), we have 
\begin{equation*}
	\begin{aligned}
		\norm{\stateSS(\initT+2\timeDown[]) - \stateTA(\initT+2\timeDown[])} & \leq \norm{\stateSS(\initT+2\timeDown[]) - \nthStateTA[1](\initT+2\timeDown[])} \\&\quad + \norm{\nthStateTA[1](\initT+2\timeDown[]) - \stateTA(\initT+2\timeDown[])} \\
		& \leq \gamBar + \constLemma\gamBar \expo[{-\lam \timeDown[]}],
	\end{aligned}
\end{equation*}
where the second inequality is due to \eqref{eq: boundIntk} and Theorem \ref{thm: thoppeResult}. 
%Similarly for \(n = 3 \), we have
%\begin{equation*}
%	\begin{aligned}
%		\norm{\stateSS(\initT+3\timeDown[]) - \stateTA(\initT+3\timeDown[])} &\leq \norm{\stateSS(\initT+3\timeDown[]) - \nthStateTA[2](\initT+3\timeDown[])} + \norm{\nthStateTA[2](3\timeDown[]) - \stateTA(3\timeDown[])} \\
%		& \leq \gamBar + \constLemma (\gamBar + \constLemma\gamBar \expo[{-\lam \timeDown[]}]) \expo[{-\lam \timeDown[]}] \\
%		& = \gamBar (1 + \constLemma\expo[{-\lam \timeDown[]}]  + \constLemma^2 \expo[{-2\lam \timeDown[]}]).
%	\end{aligned}
%\end{equation*}
Fix \(q \in \N \), and assume, as the induction hypothesis, that the claim holds for each \(n \in \setInt[q-1] \). Define \(t_q \Let \initT + q\timeDown[] \). For the induction step \(n = q \), we observe that 
\begin{equation*}
	\begin{aligned}
		\norm{\stateSS(t_q) - \stateTA(t_q)} &\leq \norm{\stateSS(t_q) - \nthStateTA[q-1](t_q)} + \norm{\nthStateTA[q-1](t_q) - \stateTA(t_q)} \\
		& \leq \gamBar + \constLemma \bigg(\gamBar \Big(\sum_{i=0}^{q-2} \constLemma^i\expo[{-i \lam \timeDown[]}]\Big) \expo[{-\lam \timeDown[]}] \bigg)
	\end{aligned}
\end{equation*}
where the second inequality is due to \eqref{eq: nthSum}, which proves the claim \eqref{eq: nthSum}. A pictorial depiction of the bounding process \eqref{eq: nthSum} is provided in Fig. \ref{fig: xplain}. Note that the right hand side of \eqref{eq: nthSum} is valid for time instants that are integral multiples of \(\timeDown[]\) starting at \(\initT \). The corresponding bound for a general time \(t \in [\initT, +\infty[ \) is given by 
\begin{equation}\label{eq: genTime}
\begin{aligned}
\norm{\stateSS(t) - \stateTA(t)}  \leq \gamBar + \constLemma \bigg(\gamBar \Big(\sum_{i=0}^{\floorFac-1} \constLemma^i\expo[{-i \lam \timeDown[]}]\Big) \bigg),
\end{aligned}
\end{equation}
where \(\floorFac \Let \flor[\frac{t-\initT}{\timeDown[]}] \) for \(t\geq \initT+\timeDown[] \), where we have employed \eqref{eq: nthSum} to arrive at \eqref{eq: genTime}. Indeed, for \(t \in \loro{\initT+\floorFac\timeDown[], \initT+(\floorFac+1)\timeDown[]} \) we have 
\begin{equation*}
	\begin{aligned}
		\norm{\stateSS(t) - \stateTA(t)} &\leq \norm{y(t) - \nthStateTA[{\floorFac}](t)} + \norm{\nthStateTA[{\floorFac}](t) - \stateTA(t)} \\
		& \leq \gamBar + \constLemma \bigg(\gamBar \Big(\sum_{i=0}^{\floorFac-1} \constLemma^i\expo[{-i \lam \timeDown[]}]\Big) \expo[{-\lam (t-\floorFac\timeDown[])}] \bigg)\\
		& \leq \gamBar + \constLemma \bigg(\gamBar \Big(\sum_{i=0}^{\floorFac-1} \constLemma^i\expo[{-i \lam \timeDown[]}]\Big) \bigg),
	\end{aligned}
\end{equation*}
where the second inequality is due to \eqref{eq: nthSum} and Theorem \ref{thm: thoppeResult}. 
\begin{figure}
	\centering
	\begin{tikzpicture}[thick,scale=0.8, every node/.style={scale=0.8}][
	thick,
	>=stealth',
	axis/.style={very thick, ->, >=stealth', line join=miter},
	important line/.style={thick}, dashed line/.style={dashed, thin},
	dot/.style = {
		draw,
		fill = white,
		circle,
		inner sep = 0pt,
		minimum size = 4pt
	}
	]
	\coordinate (O) at (0,0);
	\draw[->] (-0.3,0) -- (8.5,0) coordinate[label = {right:$t$}] (xmax);
	\draw[->] (0,-0.3) -- (0,4) coordinate (ymax);
	\draw (2, 0.1) -- (2, -0.1) node [below] {$\initT+T$};
	\foreach \x in {2,...,4}
	\draw (2*\x, 0.1) -- (2*\x, -0.1) node [below] {$\initT + \x T$};
	\draw[loosely dashed] (2,0) -- (2,2);
	\draw[loosely dashed] (4,0) -- (4,2);

	\draw[orange!80, ->] [densely dashed](0, 2) .. controls(4,2) and (4,0) .. (8,0.1);
	\draw[blue!80][densely dashed] (2,2.5) .. controls(3,2.4) .. (4,2.05);
	\draw[blue!80, ->][densely dashed] (4,2.05) .. controls(8,0) and(8,0.6) .. (8.5,0.5);
	
	\draw[red, ->] plot[smooth] coordinates {(0, 2) (2,2.5) (4,2.75) (6,2.5) (8,2.15) };
	\coordinate [label=left:$\normInit$](A) at (-0,2);
	\draw[<->] (2,2.5) -- (2,1.75) node[pos=0.5, right] {$\gamma$};
	\draw[<->] (4,2.75) -- (4,2.0) node[pos=0.5, right] {$\gamma$};
	\draw[<->] (4,2.0) -- (4,1.1) node[pos=0.8, right] {$\gamBar\expo[{-T}]$};
	\draw[orange][densely dashed] (6,4.5) -- (7,4.5) node[pos=1.15, right] {$x(t)$};
	\draw[red] (6,4) -- (7,4) node[pos=1.15, right] {$y(t)$};
	\draw[blue][densely dashed] (6,3.5) -- (7,3.5) node[pos=1.15, right] {$x_{(1)}(t)$};
	%	\draw[] (0,2.95) -- (0,2.95) node[pos=1.15, left] {$\min{(r,r^\prime)}$};
	%	
	\node[left] at (-0.05,-0.3) {$\initT$};
	\end{tikzpicture}
	\caption{Norm of solutions vs time for the case when \(K_1 = 1\)}\label{fig: xplain}
\end{figure}
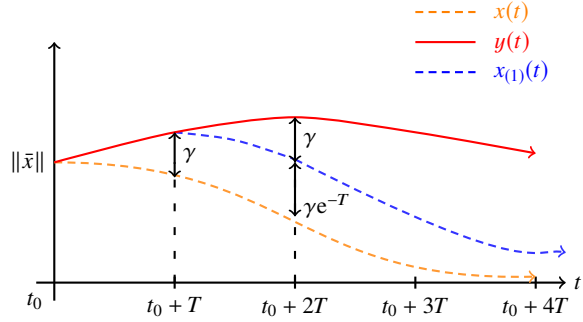

Based on the selection criterion of \(\timeDown[] \), we have \(\FAC< 1 \); thus, we can bound the expression in \eqref{eq: genTime} by the sum of the geometric series, leading to
\begin{align*}
	\bigl|\norm{\stateSS(t)} - \norm{\stateTA(t)}\bigr| \leq  \norm{\stateSS(t) - \stateTA(t)} \leq  \gamBar\Big(1 + \frac{\constLemma}{1-\FAC} \Big) \quad\text{for all}\ t \geq \initT.
\end{align*}
As a consequence, we arrive at 
\begin{align}\label{eq: StableEqn}
	\norm{\stateSS(t)} \leq \norm{\stateTA(t)} + \consLem\gamBar \leq \ups \quad\text{for all}\ t \geq \initT, 
\end{align}
where the last inequality is due to the choice of \(\gamBar \) made earlier.

This concludes the proof of Theorem \ref{thm: MainResult} since for every \(\ups> 0 \)  we have exhibited a \(\del > 0\) and a switching time \(\switchTimeS > 0\) (refer to \textbf{Step-I}), both of which are independent of the initial time \(\initT \) and depend only on \(\ups \), such that
\(\norm{\stateSS(t)} \leq \ups \quad\text{for}\ t \geq t_0 \ \text{if} \ \normInit \leq \del.  \)

\subsection{Proof of Theorem \ref{thm: AS_stable}}\label{app_label: ASstable}
This proof follows the same steps as that of Theorem \ref{thm: MainResult}, namely \ref{step:1}-\ref{step:2}.

As the equilibrium \(0 \in \R^d \) of \eqref{eq: normalSys} is asymptotically stable, we can find a \(\initAsymp > 0 \) such that if \(\normInit \leq \initAsymp \) then \(\norm{\stateTA(t)} \xrightarrow[t \rightarrow +\infty]{} 0\). Fix \(\ups > 0 \), and define \(\betBound \Let \min(\ups, \initAsymp) \). Define a compact set  \(\bddSet \Let   \ball{\betBound}{0} \). Pick a scalar \(\bound \in \loro{0,1} \) and define the another compact set \(\bddSet_\bound \Let \ball{\bound\betBound}{0}\). Also, the equilibrium point \(0 \in \R^d \) of \eqref{eq: normalSys}  is stable in the sense of Lyapunov therefore we can find a positive scalar \(\del \) such that if \(\normInit \leq \del \) then \(\stateTA(t) \in \bddSet_{\alpha} \) for all \(t \geq \initT \). In other words, the trajectory \(\stateTA(\cdot) \subset \bddSet_{\alpha} \subsetneq \bddSet \).

Pick \(\facAS \in \loro{0,1} \). For \(\timeDown[]>0 \), define \(\FAC \Let \frac{\constLemma}{\facAS} \expo[{-\lam\timeDown[]}] \), where \(K_1, \lam \) are positive scalars which depend on the set \(\bddSet\) and the closed-loop vector field governing \eqref{eq: normalSys}. Choose \(T > 0\) such that \(\FAC < 1 \). As in \S\ref{app_lab: Main Result}, we partition the interval \([\initT, +\infty[ \) into a countable family of disjoint half-open intervals of length \(\timeDown[]\) as \(\lcro{\initT, +\infty} = \bigcup_{n\in\nat} \intV[n]\) where \( \intV[n] \Let  \lcro{\initT + n \timeDown[], \initT + (n+1)\timeDown[]}\). 

As in the proof of Theorem \ref{thm: MainResult}, we appeal to Theorem \ref{thm: FilTheorem} to claim that the trajectory \(\stateAS(\cdot) \) converges uniformly to \(\stateTA(\cdot) \) on the interval \(\intV[0] \). That is, for any \(\gamBar> 0\) there exists a switching time \(\switchTime\) sufficiently small such that 
\begin{align}\label{eq: unifBoundMR2}
	\norm{\stateAS(t) - \stateTA(t)} \leq \gamBar \quad \text{for} \ t \in \intV[0].
\end{align}
At this stage pick \(\gamBar < \frac{\betBound(1 - \alpha)}{\consLem} \) with \(\consLem \Let 1 + \frac{\constLemma}{1-\FAC} \), and corresponding to this value of \(\gamBar \) we select \(\switchTimeS\) such that \eqref{eq: unifBoundMR2} is satisfied. For \(t \in \intV[0] \), we set \(\switchTimeAS(t) = \switchTimeS \). 

Analogous to the proof of Theorem \ref{thm: MainResult}, we construct an ensemble of trajectories \(\big( \nthStateTAS[k](\cdot) \big)_{k \in \N} \) each satisfying the dynamics \eqref{eq: sparseSysAS}, and the condition that  
\begin{align*}
	\nthStateTAS[n](\initT+n\timeDown[]) = \stateAS(\initT+n\timeDown[]), \quad \text{for all}\ n \in \N.
\end{align*} 

For each of the time intervals \(\intV[i] \) (\(i\in\N \)), we choose a switching time \(\switchTime_i \in \R \) such that 
\(
\norm{\stateAS(t) - \nthStateTAS[i](t)} \leq \facAS^i\gamBar \quad \text{for} \ t \in \intV[i],
\)
where a priori. Note that the existence of such \(\switchTime_i \) is guaranteed by Theorem \ref{thm: FilTheorem}. In fact, it is possible to pick \(\switchTime_i \geq \switchTime_{i+1} \) for all \(i \in \N \). We set \(\switchTimeAS(t) = \switchTime_i\) for \(t \in \intV[i]\) and \(i \in \N \). 

We claim that the difference between \(\stateAS(\cdot) \) and \(\stateTA(\cdot) \), at each time \(t\), satisfies the estimate
\begin{equation}\label{eq: genTimeAS}
\begin{aligned}
\norm{\stateAS(t) - \stateTA(t)}  \leq \facAS^{\floorFac}\gamBar + \constLemma \bigg(\gamBar \Big(\sum_{i=0}^{\floorFac-1} \facAS^{\floorFac-1-i}\constLemma^i\expo[{-i \lam \timeDown[]}]\Big) \bigg),
\end{aligned}
\end{equation}
where \(\floorFac = \flor[\frac{t-\initT}{\timeDown[]}] \) for \(t\geq \initT+\timeDown[] \). Proof of which is exactly on the lines of the proof of \eqref{eq: genTime} in \S\ref{app_lab: Main Result}. 
%For \(\floorFac \timeDown[] < t-\initT < (\floorFac+1)\timeDown[] \), we have 
%\begin{equation*}
%	\begin{aligned}
%		\norm{\stateAS(t) - \stateTA(t)} &\leq \norm{\stateAS(t) - \nthStateTAS[{\floorFac}](t)} + \norm{\nthStateTAS[{\floorFac}](t) - \stateTA(t)} \\
%		& \leq \facAS^{\floorFac}\gamBar + \constLemma \bigg(\gamBar \Big(\sum_{i=0}^{\floorFac-1} \facAS^{\floorFac-1-i}\constLemma^i\expo[{-i \lam \timeDown[]}]\Big) \expo[{-\lam (t-\floorFac\timeDown[])}] \bigg) \quad \text{\big(Using \eqref{eq: nthSum} and Theorem \ref{thm: thoppeResult} \big)}\\
%		& \leq \facAS^{\floorFac}\gamBar + \constLemma \bigg(\gamBar \Big(\sum_{i=0}^{\floorFac-1} \facAS^{\floorFac-1-i}\constLemma^i\expo[{-i \lam \timeDown[]}]\Big) \bigg).
%	\end{aligned}
%\end{equation*}

Based on the selection criterion of \(\timeDown[] \), we have \(\FAC< 1 \) and therefore we can bound  \eqref{eq: genTimeAS} as follows:
\begin{align}
	\bigl|\norm{\stateAS(t)} - \norm{\stateTA(t)}\bigr| \leq  \norm{\stateAS(t) - \stateTA(t)} \leq  \facAS^{\floorFac}\gamBar\Big(1 + \frac{\constLemma}{1-\FAC} \Big) \quad \text{for}\ t \geq \initT.
\end{align}
The above inequality suggests that 
\begin{align}\label{eq: StableEqnAS}
	\norm{\stateAS(t)} \leq \norm{\stateTA(t)} + \facAS^{\floorFac}\consLem\gamBar \quad\text{for}\ t \geq \initT. 
\end{align}

Since the equilibrium of \eqref{eq: normalSys} is locally asymptotically stable
%, there exists \(\tilde{\del}> 0\) such that for every \(\zet > 0 \) there exits \(\tilde{T} > 0 \) such that \(\|x(t)\| \leq \zet \) for all \(t \geq \initT + \tilde{T} \) given \(\normInit \leq \tilde{\del} \). Therefore 
there exists \(\del' \) such that for every \(\sca > 0 \) there exists \(T' > 0 \) such that 
\(
\normInit \leq \del' \implies \norm{x(t)} \leq \frac{\sca}{2} \ \text{for all}\ t \geq \initT+T'.
\)
Using \eqref{eq: StableEqnAS} and setting \(\normInit \leq \min({\del, \del'}) \) we can conclusively say that for every \(\sca > 0 \) if we define
\begin{align}\label{eq: TAS}
	\Tas \Let \max\bigg(T', \frac{\log(\frac{\sca}{2\consLem\gamBar})}{\log(\facAS)}\bigg), 
\end{align}
which is independent of \(\initT \), then \(\norm{\stateAS(t)} \leq \sca \) for all \(t \geq \initT + T_{AS} \). Our proof is complete.

\begin{remark}\label{rem: performanceLoss}
	{\rm It is evident from \eqref{eq: TAS} that the time required for the trajectories \(\stateAS(\cdot) \) to reach any neighborhood of the equilibrium point \(0\) is at least as large as that for the trajectories \(\stateTA(\cdot)\). This observation is quite natural since a loss of performance is expected due to the round-robin periodic scheduling of control channels.}
\end{remark} 

\begin{tcolorbox}
	For a consistent notation in Appendices \ref{ssec: NonUniform} and \ref{ssec: Uniform}, we denote the trajectories of \eqref{eq: normalSys} by \(\stateTA(\cdot;\initROA,\initT) \) and of  \eqref{eq: sparseSysAS} by \(\stateAS(\cdot;\initROA,\switchFuncAs_{\initT}(\cdot,\switchTime),\initT)\), where \(\stateTA(\initT;\initROA,\initT) = \stateAS(\initT;\initROA,\switchFuncAs_{\initT}(\cdot,\switchTime),\initT) = \initROA\) and \(\switchFuncAs_{\initT}(\cdot,\tau)\) is the switching time defined in \eqref{eq: switchRuleAS} corresponding to \(\switchTime\) where the switching starts at \(\initT\). Such an elaborate notation will be useful.
\end{tcolorbox}

{
	\subsection{Proof of Theorem \ref{thm: NonuniformResult}}
	\label{ssec: NonUniform}
	Recall from \cite[Chapter 5, Lemma 45]{vidyasagar2002nonlinear} that the region of attraction \(\roa{cl}\) of the equilibrium \(0\) of \eqref{eq: normalSys} is an open, connected and invariant set. We split the proof in two steps.
	
	\noindent\(\bullet\) \textbf{Step-I:} Pick \(c > 0\) such that \(\ball{c}{0} \subset \roa{cl} \). Such a \(c\) exists because \(\roa{cl}\) is open as noted above. If we initialize the systems \eqref{eq: normalSys} and \eqref{eq: sparseSysAS} at time \(\pseudoT > 0\) from the common initial condition \(\pseudoIC\), then Theorem \ref{thm: AS_stable} ensures existence of a pair \((\Delta_{c},\tau_c(\cdot))\) independent of \(\pseudoT\) such that for every \(\switchTime\) that dominates \(\switchTime_c\) and for every \(\pseudoIC \in \ball{{\Delta_c}}{0}\), the trajectory \(\stateAS(\cdot;\pseudoIC,\switchFuncAs_{\pseudoT}(\cdot,\switchTime),\pseudoT)\) lies entirely in \(\ball{c}{0}\) and converges to \(0\) asymptotically. If the given initial condition \(\stinit\) at \(\initT\) is in \(\ball{{\Delta_c}}{0}\), there is clearly nothing left to prove. Otherwise, we focus on the case when the initial condition \(\stinit\) at \(\initT\) does not lie inside \(\ball{\Delta_c}{0}\); this case is treated in \textbf{Step-II}.
	
	\noindent\(\bullet\) \textbf{Step-II:} We pick \(c > 0\) as in \textbf{Step-I}, select \(\alpha\in[0,1]\), and define \(\Delta_c^{(\alpha)} \defas \alpha \Delta_c\). Recall \cite[Definition 43]{vidyasagar2002nonlinear} that for every \(\stinit\) in the basin of attraction \(\roa{cl}\) of the equilibrium point \(0\) of \eqref{eq: normalSys}, there exists \(T(\stinit) > 0\) such that for all \(t\geq \initT+T(\stinit)\) we have \(\stateTA(t;\stinit,\initT) \in  \ball{{\Delta_c^{(\alpha)}}}{0}\). Theorem \ref{thm: FilTheorem} asserts that for \(\gamma = \Delta_c(1-\alpha)\) there exists \(\pseudoST > 0\) (sufficiently small and constant but depending on \(\stinit, \Delta_c, \alpha\)) such that 
	\begin{align}\label{eq: NonUniform}
		\|\stateAS(t;\stinit,\switchFuncAs_{\initT}(\cdot,\pseudoST),\initT)-\stateTA(t;\stinit,\initT)\| \leq \gamma \quad \text{for} \ t \in \lcro{\initT,\initT+T(\bar{x})}.
	\end{align}
	%	With slight abuse of notation, on the interval \(\lcro{\initT,\initT+T(\bar{x})}\) the switch time is taken to be constant \(\pseudoST\). 
	The triangle inequality and continuity of solutions shows us, in particular, that 
	\[
	z(\initT+T(\bar{x});\stinit,\switchFuncAs_{\initT}(\cdot,\pseudoST),\initT) \in \ball{{\Delta_c}}{0}.
	\]
	In other words, once the trajectory of \eqref{eq: sparseSysAS} initialized at \(\stinit\) and scheduled in accordance with \(\switchFuncAs_{\initT}(\cdot,\pseudoST)\) reaches \(\ball{{\Delta_c}}{0}\), we transition to the switching discussed in \textbf{Step-I}. 
	
	In summary, the convergence of \eqref{eq: sparseSysAS} to \(0\) is achieved by employing any switching time that dominates the following function:%first switch with constant \(\pseudoST\) (infact we can switch faster than this as well but to avoid more notations we avoid stating that repeatedly) during \([\initT,\initT+T(\bar{x})[\) and then following it up with \(\switchTime_c(\cdot)\). 
	% Recall that \(\switchTime_c\) is the switching signal which was discussed in \textbf{(Step-I)}. 
	%Thus, if we schedule \eqref{eq: sparseSysAS} with any switching time \(\switchTimeAS\) that dominates
	\[
	\switchTime'(\cdot;\initROA)\defas \begin{cases}
	\pseudoST \quad  &\text{if}\ t \in \lcro{\initT,\initT+T(\initROA)},\\
	\switchTime_c(\cdot) \quad &\text{if} \ t \geq \initT+T(\bar{x}),
	\end{cases}
	\]
	where we set \(T(\initROA) = 0\) whenever \(\stinit\in\ball{{\Delta_c}}{0}\).
	
	\subsection{Proof of Theorem \ref{thm: RegionOfAttractionUniformResult}}\label{ssec: Uniform}
	We follow the same route as in Appendix \ref{ssec: NonUniform}.
	
	\noindent\(\bullet\) \textbf{Step-I:} Let \(\rbar \defas \sup\bigl\{r > 0\,\big|\,\ball{r}{0}\subset \roa{cl}\bigr\}\). If we initialize the systems \eqref{eq: normalSys} and \eqref{eq: sparseSysAS} at time \(\pseudoT\) from the common initial condition \(\pseudoIC\), then Theorem \ref{thm: AS_stable} ensures the existence of a pair \((\Delta_{\rbar},\switchTime_{\rbar}(\cdot))\) such that for every \(\switchTime\) that dominates \(\switchTime_{\rbar} \) and for every \( \pseudoIC\in \ball{\Delta_{\rbar}}{0}\), we have \(z(\cdot;\pseudoIC,\switchFuncAs_{\pseudoT}(\cdot,\switchTime),\pseudoT) \subset \interior(\ball{{\rbar}}{0}) \) and \(z(\cdot;\pseudoIC,\switchFuncAs_{\pseudoT}(\cdot,\switchTime),\pseudoT)\) converges to the equilibrium \(0\) asymptotically. Thus, referring to the notation established before Theorem \ref{thm: RegionOfAttractionUniformResult}, if \(\bar{\radCompact} >0\) is such that \(\compactSet{\bar{\radCompact}} \subset \interior(\ball{{\Delta_{\rbar}}}{0}) \) and we pick \(\pseudoIC\in \compactSet{\bar{\radCompact}} \) as the initial condition for \eqref{eq: sparseSysAS}, the trajectories converge to the equilibrium with the switching time \(\switchFuncAs_{\pseudoT}(\cdot,\switchTime) \) for all \(\switchTime\) that dominates \(\switchTime_{\rbar} \). 
	
	\noindent\(\bullet\) \textbf{Step-II:} We consider the case of \(\radCompact > 0\) such that \(\compactSet{\radCompact}\not\subset\interior(\ball{{\Delta_{\rbar}}}{0}) \). For such \(\radCompact\) we pick \( \tilde{\Delta} < \Delta_{\rbar} \) such that \(\ball{{\tilde{\Delta}}}{0} \subset \compactSet{\radCompact}\). Next pick \(\alpha \in \loro{0,1}\) and fix  \(\tilde{\radCompact} < \radCompact\) such that \(\compactSet{\tilde{\theta}} \subset \ball{{\alpha\tilde{\Delta}}}{0}\). Define
	%Then for \(\initROA \in \compactSet{\radCompact} \)  and \(\initROA\not\in \compactSet{\tilde{\theta}}\) define 
	\[
	\tilde{T} \defas \frac{\radCompact-\tilde{\radCompact}}{\inf_{x\in \compactSet{\radCompact}\setminus\compactSet{\tilde{\radCompact}}}|\inprod{\nabla \lyap(x)}{F(x)}|}
	\]
	where \(F(x) \Let \drift[x] + \sumFunc{k}{\contDim} \channel{k}{x} \control[k](x) \) is the RHS of \eqref{eq: normalSys} and \(\lyap\) is the Lyapunov function for \eqref{eq: normalSys} described before Theorem \ref{thm: RegionOfAttractionUniformResult}. Note that the denominator in the definition of \(\tilde{T}\) is strictly positive by definition of \(V\), and therefore \(\tilde{T} \) is well defined. We note that \(\tilde{T}\) depends only on \(\radCompact\) in addition to system dependent objects, and we shall denote this dependence on \(\radCompact\) explicitly in the sequel. Thus, for any \(\stinit \in \compactSet{\theta}\backslash \compactSet{\tilde{\theta}} \), we have \(\stateTA(t;\initROA,\initT) \in \ball{{\compactSet{\tilde{\theta}}}}{0} \) for all \(t\geq \initT + \tilde{T}(\radCompact) \). Theorem \ref{thm: FilTheorem} asserts that for \(\gamma = \tilde{\Delta}(1-\alpha) \) there exists \(\pseudoST > 0\) (constant and sufficiently small) such that 
	\begin{align}\label{eq: Uniform}
		\|\stateAS(t;\stinit,\switchFuncAs_{\initT}(\cdot,\pseudoST),\initT) - \stateTA(t;\stinit,\initT)\| \leq \gamma \quad \text{for} \ t \in \lcro{\initT,\initT+\tilde{T}(\radCompact)}.
	\end{align}
	%With slight abuse of notation, on the interval \(\lcro{\initT,\initT+T(\bar{x})}\) the switch time is taken to be constant \(\pseudoST\). 
	The triangle inequality, the estimate \eqref{eq: Uniform}, and continuity of trajectories shows that
	\[
	\stateAS(\initT+\tilde{T}(\radCompact);\stinit,\switchFuncAs_{\initT}(\cdot,\pseudoST);\initT ) \in \ball{{\tilde{\Delta}}}{0}.
	\]
	Once the scheduled trajectory reaches \(\ball{{\tilde{\Delta}}}{0}\), we follow the switching time discussed in \textbf{Step-I}.
	
	In summary, the convergence of \eqref{eq: sparseSysAS} to the equilibrium \(0\) is achieved by using any switching time \(\switchTimeAS(\cdot)\) that dominates the following function:
	\[
	\switchTime'(\cdot;\radCompact)\defas
	\begin{cases}
	\pseudoST \quad  &\text{if}\ t \in \lcro{\initT,\initT+\tilde{T}(\radCompact)},\\
	\switchTime_{\rbar}(\cdot) \quad &\text{if} \ t \geq \initT+\tilde{T}(\radCompact),
	\end{cases}
	\]
	where we set \(\tilde{T}(\radCompact) = 0\) whenever \(\compactSet{\radCompact}\subset\interior(\ball{{\Delta_{\rbar}}}{0})\).
	%then the resulting trajectory converges to the equilibrium.
}

\bibliography{references}

\newcommand{\etalchar}[1]{$^{#1}$}
\begin{thebibliography}{WPNH19}

\bibitem[Art08]{ref:Art-08}
Z.~Artstein.
\newblock A {Y}oung measures approach to averaging.
\newblock In {\em Differential Equations, Chaos and Variational Problems},
  volume~75 of {\em Progress in Nonlinear Differential Equations and their
  Applications}, pages 15--28. Birkh\"auser, Basel, 2008.

\bibitem[BM10]{ref:BacMaz-10}
A.~Bacciotti and L.~Mazzi.
\newblock Stabilisability of nonlinear systems by means of time-dependent
  switching rules.
\newblock {\em International Journal of Control}, 83(4):810--815, 2010.

\bibitem[BM13]{ref:BerMed-13}
L.~D. Berkovitz and N.~G. Medhin.
\newblock {\em Nonlinear {O}ptimal {C}ontrol {T}heory}.
\newblock Chapman \& Hall/CRC Applied Mathematics and Nonlinear Science Series.
  CRC Press, Boca Raton, FL, 2013.

\bibitem[CC10]{chen2010review}
B.~Chen and H.~Cheng.
\newblock A review of the applications of agent technology in traffic and
  transportation systems.
\newblock {\em IEEE Transactions on Intelligent Transportation Systems},
  11(2):485--497, 2010.

\bibitem[CHW11]{claes2011decentralized}
R.~Claes, T.~Holvoet, and D.~Weyns.
\newblock A decentralized approach for anticipatory vehicle routing using
  delegate multi agent systems.
\newblock {\em IEEE Transactions on Intelligent Transportation Systems},
  12(2):364--373, 2011.

\bibitem[CM12]{SG1}
A.~Chakrabortty and D.~Marija, editors.
\newblock {\em Control and {O}ptimization {M}ethods for {E}lectric {S}mart
  {G}rids}.
\newblock Power Electronics and Power Systems. Springer, New York, 2012.

\bibitem[CPRT17]{JurdQuinn}
M.~Caponigro, B.~Piccoli, F.~Rossi, and E.~Tr\'{e}lat.
\newblock Sparse {J}urdjevic-{Q}uinn stabilization of dissipative systems.
\newblock {\em Automatica}, 86:110--120, 2017.

\bibitem[DHVH11]{donkers2011stability}
M.~C.~F. Donkers, W.~M. P.~H. Heemels, N.~{Van de Wouw}, and L.~Hetel.
\newblock Stability analysis of networked control systems using a switched
  linear systems approach.
\newblock {\em IEEE {T}ransactions on {A}utomatic control}, 56(9):2101--2115,
  2011.

\bibitem[Fil88]{filippov2013differential}
A.~F. Filippov.
\newblock {\em Differential {E}quations with {D}iscontinuous {R}ighthand
  {S}ides}, volume~18 of {\em Mathematics and its Applications (Soviet
  Series)}.
\newblock Kluwer Academic Publishers Group, Dordrecht, 1988.
\newblock Translated from the Russian.

\bibitem[GB08]{el2008chaotic}
F.~E. Guezar and H.~Bouzahir.
\newblock Chaotic behavior in a switched dynamical system.
\newblock {\em Modelling and Simulation in Engineering}, 2008:2, 2008.

\bibitem[GIL07]{gorges2007optimal}
D.~G\"{o}rges, M.~Iz\'ak, and D.~Liu.
\newblock Optimal control of systems with resource constraints.
\newblock In {\em Proceedings of the 46th IEEE Conference on Decision and
  Control}, pages 1070--1075, 2007.

\bibitem[IYB06]{YukselOptimalControl}
O.~C. Imer, S.~Y\"{u}ksel, and T.~Ba\c{s}ar.
\newblock Optimal control of {LTI} systems over unreliable communication links.
\newblock {\em Automatica}, 42(9):1429--1439, 2006.

\bibitem[JQ78]{jurdjevic1978controllability}
V.~Jurdjevic and J.~P. Quinn.
\newblock Controllability and stability.
\newblock {\em Journal of {D}ifferential {E}quations}, 28(3):381--389, 1978.

\bibitem[KC17]{ref:KunCha-17}
A.~Kundu and D.~Chatterjee.
\newblock Stabilizing switching signals: a transition from point-wise to
  asymptotic conditions.
\newblock {\em Systems \& Control Letters}, 106:16--23, 2017.

\bibitem[LA05]{ref:astolfi-02}
M.~Lovera and A.~Astolfi.
\newblock Global magnetic attitude control of inertially pointing spacecraft.
\newblock {\em Journal of {G}uidance, {C}ontrol, and {D}ynamics},
  28(5):1065--1072, 2005.

\bibitem[Lib03]{LibBook}
D.~Liberzon.
\newblock {\em Switching in {S}ystems and {C}ontrol}.
\newblock Systems \& Control: Foundations \& Applications. Birkh\"{a}user
  Boston, Inc., Boston, MA, 2003.

\bibitem[MSBH17]{van2017switched}
A.~V. Maas, Y.~F. Steinbuch, A.~Boverhof, and W.~P. Heemels.
\newblock Switched control of a {SCARA} robot with shared actuation resources.
\newblock {\em IFAC-PapersOnLine}, 50(1):1931--1936, 2017.

\bibitem[NT04a]{Nesic2004TAC}
D.~Ne\v{s}i\'{c} and A.~R. Teel.
\newblock Input-output stability properties of networked control systems.
\newblock {\em IEEE Trans. Automat. Control}, 49(10):1650--1667, 2004.

\bibitem[NT04b]{Nesic2004Automatica}
D.~Ne\v{s}i\'{c} and A.~R. Teel.
\newblock Input-to-state stability of networked control systems.
\newblock {\em Automatica J. IFAC}, 40(12):2121--2128 (2005), 2004.

\bibitem[PKS13]{polyak2013lmi}
B.~Polyak, M.~Khlebnikov, and P.~Shcherbakov.
\newblock An {LMI} approach to structured sparse feedback design in linear
  control systems.
\newblock In {\em Proceedings of the European Control Conference}, pages
  833--838, 2013.

\bibitem[PL01]{panteley2001uniform}
E.~Panteley and A.~Lor{\'\i}a.
\newblock Uniform exponential stability for parameterized linear
  ``skew-symmetric'' systems.
\newblock In {\em Proceedings of European Control Conference}, pages
  2410--2415. IEEE, 2001.

\bibitem[Sch09]{SchenatoZeroing}
L.~Schenato.
\newblock To zero or to hold control inputs with lossy links?
\newblock {\em IEEE Transactions on Automatic Control}, 54(5):1093--1099, 2009.

\bibitem[SSF{\etalchar{+}}07]{schenato2007foundations}
L.~Schenato, B.~Sinopoli, M.~Franceschetti, K.~Poolla, and S.~S. Sastry.
\newblock Foundations of control and estimation over lossy networks.
\newblock {\em Proceedings of the IEEE}, 95(1):163--187, 2007.

\bibitem[TB19]{ConcentrationBound}
G.~Thoppe and V.~Borkar.
\newblock A concentration bound for stochastic approximation via {A}lekseev's
  formula.
\newblock {\em Stochastic Systems}, 9(1):1--26, 2019.

\bibitem[Tes12]{teschl2012ordinary}
G.~Teschl.
\newblock {\em {O}rdinary {D}ifferential {E}quations and {D}ynamical
  {S}ystems}, volume 140.
\newblock American Mathematical Society, 2012.

\bibitem[Vid02]{vidyasagar2002nonlinear}
M.~Vidyasagar.
\newblock {\em Nonlinear {S}ystems {A}nalysis}.
\newblock Classics in {A}pplied {M}athematics. SIAM, 2nd edition, 2002.

\bibitem[WBB01]{walsh2001asymptotic}
G.~C. Walsh, O.~Beldiman, and L.~G. Bushnell.
\newblock Asymptotic behavior of nonlinear networked control systems.
\newblock {\em IEEE {T}ransactions on {A}utomatic {C}ontrol}, 46(7):1093--1097,
  2001.

\bibitem[WD91]{ref:wen-99}
J.~T. Wen and K.~K. Delgado.
\newblock The attitude control problem.
\newblock {\em IEEE Transactions on Automatic control}, 36(10):1148--1162,
  1991.

\bibitem[WNP15]{wang2015emulation}
W.~Wang, D.~Ne{\v{s}}i{\'c}, and R.~Postoyan.
\newblock Emulation-based stabilization of networked control systems
  implemented on {F}lex{R}ay.
\newblock {\em Automatica}, 59:73--83, 2015.

\bibitem[WNT12]{wang2012input}
W.~Wang, D.~Ne{\v{s}}i{\'c}, and A.~R. Teel.
\newblock Input-to-state stability for a class of hybrid dynamical systems via
  averaging.
\newblock {\em {M}athematics of {C}ontrol, {S}ignals, and {S}ystems},
  23(4):223--256, 2012.

\bibitem[WPNH19]{Heemels2019periodic}
W.~Wang, R.~Postoyan, D.~Ne\v{s}i\'{c}, and W.~Heemels.
\newblock Periodic event-triggered control for nonlinear networked control
  systems.
\newblock {\em IEEE {T}ransactions on {A}utomatic {C}ontrol}, 2019.

\bibitem[WY01]{walsh2001scheduling}
G.~C. Walsh and H.~Ye.
\newblock Scheduling of networked control systems.
\newblock {\em IEEE {C}ontrol {S}ystems {M}agazine}, 21(1):57--65, 2001.

\end{thebibliography}
\bibliographystyle{alpha}
\end{document}